\documentclass{article}
\usepackage{amsmath,amssymb}
\usepackage{mathptmx}
\usepackage{graphicx}
\usepackage{algorithm,algorithmic}
\usepackage{placeins}
\usepackage[labelformat=simple]{subfig}

\title{Moving Mesh Discontinuous Galerkin Methods for PDEs with Traveling Waves}
\author{ Murat Uzunca, B\"{u}lent Karas\"{o}zen, Tu\u{g}ba K\"u\c{c}\"ukseyhan}
\date{}

\begin{document}

\maketitle

\begin{abstract}
In this paper, a moving mesh discontinuous Galerkin (dG) method is developed for nonlinear partial differential equations (PDEs) with traveling wave solutions. The moving mesh strategy for one dimensional PDEs is based on the rezoning approach which decouples the solution of the PDE from the moving mesh equation. We show that the dG moving mesh method is able to resolve sharp wave fronts and wave speeds accurately for the optimal, arc-length and curvature monitor functions.  Numerical results reveal the efficiency of the proposed moving mesh dG method for solving Burgers', Burgers'-Fisher and Schl\"ogl(Nagumo) equations.
\end{abstract}

\section{Introduction}
\label{sec:intro}

The discontinuous Galerkin (dG) method is one of the most powerful discretization techniques for solving partial differential equations (PDEs) \cite{arnold82ipf,riviere08dgm}, especially for convection dominated problems, exhibiting localized phenomena like  sharp traveling wave fronts, internal and boundary layers \cite{karasozen14tsa,uzunca14adg}. The dG method has been applied for this kind of singularly perturbed linear and nonlinear PDEs  extensively using  h-adaptive (refinement and coarsening in space), p-adaptive (enrichment of the local polynomial degree), hp-adaptive and space-time adaptive methods in the last two decades.  Another approach to deal with these kind of problems, is the r-method or moving mesh method. In the moving mesh method the grid points are relocated in the regions where the solution shows rapid variation, while keeping the number of the nodes fixed. The dG discretization is very flexible, since there is no continuity requirement between the inter-element boundaries, which makes it  suitable as a  moving mesh method on  irregular meshes. Most of the studies with moving mesh methods are limited to  finite difference and continuous finite element discretization \cite{Huang11}. There are only few publications dealing with dG moving mesh method. They include the interior penalty dG method for preprocessing the solutions  of steady state diffusion-convection-reaction equations  \cite{Antonietti08}, and the  local dG moving mesh method for hyperbolic conservation laws \cite{Ruo06}.

In this paper we develop  an adaptive  dG moving mesh method for one dimensional semi-linear differential equations with traveling  wave solutions of the form
\begin{subequations}\label{model}
\begin{align}
  u_t  & =   \epsilon u_{xx}   - f(u, u_x), & (x,t)\in\Omega\times (t_0,T_f] \\
   u(x_L,t) &=  u_L ,  \; u(x_R,t) = u_R,  & t\in (t_0,T_f]\\
	u(x,t_0) &= u_0 , & x\in\Omega,
\end{align}
\end{subequations}
where $\Omega =[x_L,x_R]\subset\mathbb{R}$, $t_0$ and $T_f$ are the initial and final time instances, respectively, and $\epsilon$ denotes the diffusion coefficient. The model equation \eqref{model} becomes the Burgers equation with $f(u,u_x) = uu_x$ \cite{Soheili12anm} , Burgers'-Fisher equation with $f(u,u_x) = \alpha uu_x + \beta u(u-1)$  \cite{Soheili12anm} and the Schl\"ogl or Nagumo equation with $f(u,u_x) = u(1-u)(1-\beta)/\delta$ \cite{Huang11}.

A moving mesh method has three main components; the discretization of the physical PDE,  mesh strategy using monitor functions and discretization of the mesh equation. The discretization of physical PDE is either coupled with the moving mesh equation or separated. In the quasi-Lagrangian approach, a large system of the discretized PDE and moving mesh equation are solved simultaneously by the standard ordinary differential equation (ODE) solvers. Instead, we use the rezoning approach by solving alternately the PDE and mesh equation, which allows more flexibility; mesh generation can be coded separately and embedded in the solution of the PDE. Since the mesh is updated at each time step, the physical PDE has to be discretized at the next time step on the new mesh. We use the static rezoning approach with the same number of points at each time step \cite{Hyman89} in contrast to the dynamic rezoning \cite{Hyman03} where the number of mesh points is changed at every time step. Therefore in the static rezoning approach the solutions from old to new mesh have to be interpolated.

The paper is organized as follows. In the next section we describe briefly  the dG method for the 1D model problem \eqref{model} on a uniform fixed mesh. Moving mesh adaption strategy and the adaptive moving mesh dG algorithm is presented in Section~\ref{sec:moving}. Numerical results are given in Section~\ref{sec:numeric} to demonstrate the effectiveness of the proposed method.

\section{Discretization of the problem on a fixed mesh}
\label{sec:dg}

Before giving the moving mesh strategy in Section~\ref{sec:moving}, in this section we describe the dG discretization of the model problem \eqref{model} on a fixed uniform mesh
\begin{equation}\label{fix_mesh}
\mathcal{T}_h : \quad x_n=x_L+nh, \quad n=0,1,\ldots , N_I,
\end{equation}
consisting of $N_I$ elements (sub-intervals) $I_n=[x_{n-1},x_n]$, $n=1,2,\ldots , N_I$, and with the fixed mesh size $h=(x_R-x_L)/N_I$.

\subsection{Space discretization by discontinuous Galerkin method}

\begin{figure}[htb!]
\centering
\includegraphics[width=0.7\textwidth]{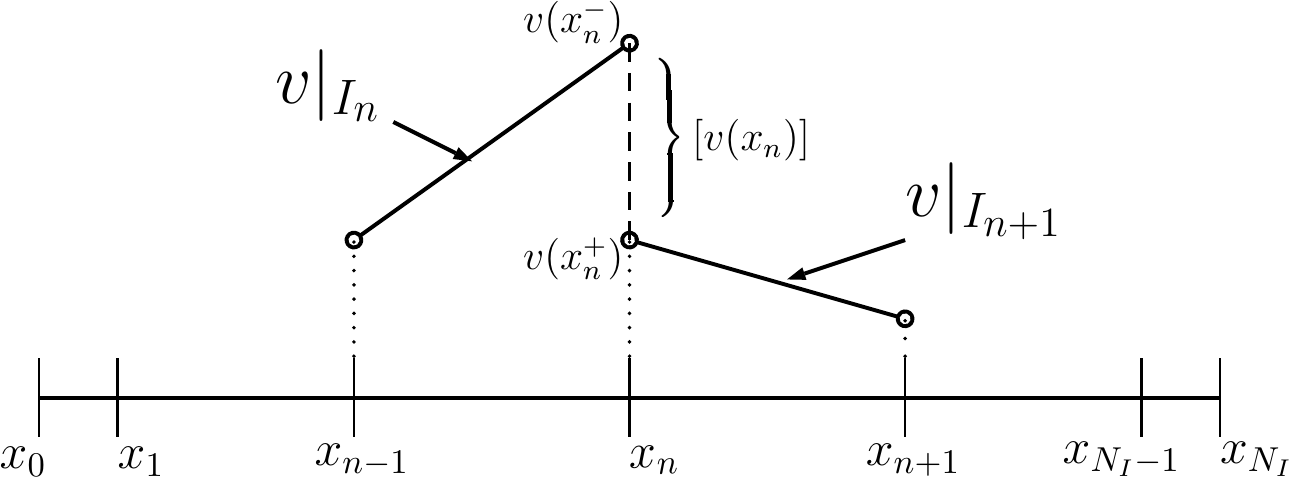}
\caption{Jump and limit terms of a piecewise discontinuous function $v(x)$.\label{fig_jump}}
\end{figure}

We use for the space discretization of the model problem \eqref{model} on a fixed mesh \eqref{fix_mesh} the symmetric interior penalty Galerkin (SIPG) method  \cite{arnold82ipf,riviere08dgm} which is a member of the family of dG methods. The dG methods use the space of piecewise discontinuous polynomials of degree at most $k$:
$$
V_h = \{ v: v|_{I_n}\in\mathbb{P}_k(I_n) \; , \; \forall n=1,\ldots , N_I\},
$$
where $\mathbb{P}_k(I_n)$ is the space of polynomials of degree at most $k$ on an interval $I_n$. Since the functions in $V_h$ are discontinuous at the inter-element nodes, we define the jump and average of a piecewise function $v$ at the endpoints of $I_n$, $n=1,\ldots , N_I-1$, respectively, as depicted in Figure~\ref{fig_jump},
\begin{equation}\label{jump1}
[v(x_n)]=v(x_n^-) - v(x_n^+) \; , \quad \{ v(x_n) \}=\frac{1}{2}( v(x_n^-) + v(x_n^+) ),
\end{equation}
with
\begin{equation}\label{limit}
v(x_n^-) = \lim_{x\mapsto x_n^-} v(x) \; , \quad v(x_n^+) = \lim_{x\mapsto x_n^+} v(x).
\end{equation}
On the boundary nodes, the jump and average are defined as
\begin{equation}\label{jump2}
[v(x_0)]=-v(x_0^+), \; \{ v(x_0) \}= v(x_0^+), \; [v(x_{N_I})]=v(x_{N_I}^-), \; \{ v(x_{N_I}) \}= v(x_{N_I}^+).
\end{equation}

The SIPG scheme is constructed by multiplying  the continuous (the solution $u$ is sufficiently smooth at the end points of $I_n$) equation \eqref{model} by a test function $v\in V_h$ and integrating by parts on each element $I_n$, $n=1,\ldots , N_I$:
$$
\int_{x_{n-1}}^{x_n}u_tvdx + \int_{x_{n-1}}^{x_n}\epsilon u_xv_x dx -\epsilon u_x(x_{n})v(x_n^-) + \epsilon u_x(x_{n-1})v(x_{n-1}^+) + \int_{x_{n-1}}^{x_n} f(u,u_x)vdx = 0.
$$
By adding all $N_I$ equations, and using the definition of the jumps  \eqref{jump1} and \eqref{jump2}, we obtain
$$
\sum_{n=1}^{N_I} \left( \int_{x_{n-1}}^{x_n}u_tvdx + \int_{x_{n-1}}^{x_n}\epsilon u_xv_x dx + \int_{x_{n-1}}^{x_n} f(u,u_x)vdx \right) - \sum_{n=0}^{N_I} [\epsilon u_x(x_{n})v(x_n)]  = 0.
$$
One can verify that for $1\leq n \leq {N_I}-1$
\begin{equation}\label{jump_id}
[\epsilon u_x(x_{n})v(x_n)] = \{\epsilon u_x(x_{n})\}[v(x_n)] + [\epsilon u_x(x_{n})]\{ v(x_n)\}.
\end{equation}
Using the identity \eqref{jump_id} and the fact that $[\epsilon u_x(x_{n})]=0$ for all $1\leq n \leq {N_I}-1$ ($u$ was sufficiently smooth at the end points of $I_n$), we obtain
\begin{equation}\label{predg}
\sum_{n=1}^{N_I} \left( \int_{x_{n-1}}^{x_n}u_tvdx + \int_{x_{n-1}}^{x_n}\epsilon u_xv_x dx + \int_{x_{n-1}}^{x_n} f(u,u_x)vdx \right) - \sum_{n=0}^{N_I} \{\epsilon u_x(x_{n})\}[v(x_n)] = 0.
\end{equation}
Additionally, we have $[u(x_n)]=0$ for all $1\leq n \leq {N_I}-1$. Then, adding the penalizing terms and the terms on the boundary nodes, $n=\{0,{N_I}\}$, to  both sides of \eqref{predg} by keeping them unknown on the left hand side and imposing the boundary conditions $u_L$ and $u_R$ on the right hand side, leads to the SIPG formulation:
\begin{equation}\label{sipg}
\int_{x_L}^{x_R} u_tvdx + a(u,v) + \int_{x_L}^{x_R} f(u,u_x)vdx = l(v).
\end{equation}
In \eqref{sipg}, $a(u,v)$ and $l(v)$ denote the symmetric bilinear form  and the linear right hand side of the SIPG scheme
\begin{align*}
a(u,v) &= \int_{x_L}^{x_R} \epsilon u_xv_x dx + \sum_{n=0}^{N_I}\left(  - \{\epsilon u_x(x_{n})\}[v(x_n)] - \{\epsilon v_x(x_{n})\}[u(x_n)] +  \frac{\sigma}{h}[u(x_{n})][v(x_n)] \right), \\
l_h(v_h) &= u_L \left( \epsilon v_x(x_0) - \frac{\sigma}{h}v(x_0) \right) + u_R \left( \frac{\sigma}{h}v(x_{N_I}) - \epsilon v_x(x_{N_I}) \right),
\end{align*}
where $\sigma$ is the penalty parameter which should be sufficiently large \cite{riviere08dgm} in order to ensure the coercivity of the bilinear form. Hence, SIPG semi-discrete form  of \eqref{model} reads as: a.e. $t\in (t_0,T_f]$, for all $v_h\in V_h$, find $u_h:=u_h(x,t)\in V_h$ such that
\begin{subequations}
\begin{align}
\int_{x_L}^{x_R} u_h(x,t_0)v_h dx &= \int_{x_L}^{x_R} u_0 v_h dx,  \label{dg0} \\
\int_{x_L}^{x_R} (\partial_t u_h)v_hdx + a(u_h,v_h) + \int_{x_L}^{x_R} f(u_h,u_{h,x})v_hdx &= l(v_h) .\label{dg}
\end{align}
\end{subequations}

\subsection{Full discretization}

Let $\{\phi_i^n\}$ denote the basis functions spanning the dG finite elements space $V_h$ for $i=1,\ldots , N_k$ and $n=1,\ldots , N_I$, where $N_k$ stands for the local dimension depending on the polynomial order $k$ and is given by $N_k=k+1$ in 1D. The local nature of dG methods leads to the basis functions and the approximate solution of the form
\begin{equation}\label{comb}
\phi_i^n(x)=\left\{
\begin{array}{cl}
\phi_i^n(x) \; , & x\in I_n \\
0 \; , & x\notin I_n
\end{array}\right. \; , \qquad u_h(t)=\sum_{n=1}^{N_I}\sum_{i=1}^{N_k} \upsilon_i^n(t)\phi_i^n(x),
\end{equation}
where $\{\upsilon_i^n (t)\}$ are the time-dependent unknown coefficients. We substitute the second identity in \eqref{comb} into the system \eqref{dg} and we choose $v_h=\phi_i^n$ for $i=1,\ldots , N_k$ and $n=1,\ldots , N_I$,  which leads to the $N=N_k\times N_I$ dimensional non-linear system of equations of \eqref{dg} in matrix-vector form
\begin{equation}\label{matrix}
M\upsilon_t + S\upsilon + h(\upsilon) - d  = 0,
\end{equation}
where $M\in\mathbb{R}^{N\times N}$ is the usual mass matrix and $S\in\mathbb{R}^{N\times N}$ is the stiffness matrix related to the bilinear form $a(u_h,v_h)$. The vectors $h(\upsilon)\in\mathbb{R}^{N}$ and $d\in\mathbb{R}^{N}$ are the non-linear vector of unknown coefficients $\upsilon$ corresponding to the integral of the nonlinear term in \eqref{dg} and the right hand side vector related to the linear form $l(v_h)$, respectively. The initial vector $\upsilon (0)$ is found by using the equation \eqref{dg0} and the second identity in \eqref{comb}.

We solve the fully discrete system of \eqref{model}, by applying the  backward Euler scheme to the semi-discrete system \eqref{matrix}. Let $t_0 < t_1 < \ldots < t_J=T_f$ be the uniform partition of the time interval $[t_0,T_f]$ into J time-steps $(t_{j-1},t_j]$, $j=1,\ldots ,J$, with the step size $\Delta t = (T_f-t_0)/J$. Let us denote the approximate coefficient vector $\upsilon (t)$ of \eqref{matrix} at the time $t=t_j$ by $\upsilon^j$. Then, the fully discrete formulation of the model \eqref{model} is given as: for all $j=1,\ldots , J$, find $\upsilon^j\in\mathbb{R}^N$ such that
\begin{eqnarray}\label{time}
M\left( \frac{\upsilon^j-\upsilon^{j-1}}{\Delta t}\right) + S\upsilon^j + h(\upsilon^j)-d^j=0,
\end{eqnarray}
which is solved by Newton's method.

\section{Adaptive moving mesh method}
\label{sec:moving}

In an adaptive moving mesh method the spatial mesh $\mathcal{T}_h$ in \eqref{fix_mesh} varies with time
\begin{equation}\label{mov_mesh}
\mathcal{T}_h(t) : \quad x_n=x(\xi_n,t), \quad n=0,1,\ldots , N_I,
\end{equation}
consisting of $N_I$ elements $I_n=[x_{n-1},x_n]$, $n=1,2,\ldots , N_I$, of  the mesh size $h_n=x_n-x_{n-1}$. In \eqref{mov_mesh} the nodes $\xi_n$ belong to the fixed and uniform mesh
\begin{equation}\label{fixmov_mesh}
\mathcal{T}_h^c : \quad \xi_n=\frac{n}{N_I}, \quad n=0,1,\ldots , N_I
\end{equation}
on the auxiliary unit interval $\Omega_c=[0,1]$, together with the boundary conditions $x(0,t)=x_L$ and $x(1,t)=x_R$. Thus, in the  adaptive moving mesh method, the mesh is adjusted as the time progresses in such a way that the mesh sizes $h_n$ get smaller in the sub-intervals where the error is large, while the mesh sizes $h_n$ are made coarser in the remaining part of the interval. The error indicator is chosen to relocate the mesh points where the solution shows large variations,  based  on  the equidistribution principle, where a mesh $\mathcal{T}_h(t)$ with the mesh points $x_L=x_0<x_1<\ldots <x_{N_I-1}<x_{N_I}=x_R$ is determined by satisfying the following relation
\begin{equation}\label{equi}
\int\limits_{x_0}^{x_1} \rho(x,t)\mathrm{d}x=\cdots =\int\limits_{x_{N_I-1}}^{x_{N_I}} \rho(x,t)\mathrm{d}x.
\end{equation}
 In this way  a continuous function $\rho(x,t) > 0$ on the interval $[x_L,x_R]$ can be distributed among the sub-intervals $I_n=[x_{n-1},x_n]$, $n=1,\ldots ,N_I$. In \eqref{equi}, the function $\rho(x,t)$ is called the mesh density function, or monitor function, choice of which stands as the key point for an adaptive moving mesh method. The most popular choices for the monitor functions $\rho:=\rho(x,t)$ are \cite{Huang11amm}
\begin{itemize}
	\item {\bf optimal}
	\begin{equation}\label{density_opt}
		\rho = \left( 1+\frac{1}{\alpha}\left |u_{xx}\right |^2\right)^{1/3},
	\end{equation}
	
	\item {\bf arc-length}
	\begin{equation}\label{density_arc}
		\rho = \left ( 1 + \left | u_x \right |^2    \right )^{1/2},
	\end{equation}
	
	\item {\bf curvature}
	\begin{equation}\label{density_curv}
		\rho = \left ( 1 + \left | u_{xx} \right|^2 \right )^{1/4},
	\end{equation}
\end{itemize}
with the intensity parameter
\begin{equation}\label{density_alpha}
\alpha =\max \left\{ 1,\left( \frac{1}{x_R-x_L}\int_{x_L}^{x_R}\left |u_{xx}\right |^{2/3}dx\right)^{3}\right\}.
\end{equation}

Finding a proper mesh $\mathcal{T}_h(t)$ using  the equidistribution condition \eqref{equi} results in a system of so-called moving mesh partial differential equation (MMPDE) \cite{huang94,Huang11amm}
\begin{subequations}\label{meshdensfunc}
\begin{align}
\frac{\partial x}{\partial t} &= \frac{1}{\tau \rho}  \frac{\partial}{\partial \xi} \left(\rho \frac{\partial x}{\partial \xi}  \right), & (\xi ,t)\in\Omega_c\times (t_0,T_f],\\
x(0,t) &= x_L,\quad x(1,t)=x_R, & t\in (t_0,T_f].
\end{align}
\end{subequations}
The system \eqref{meshdensfunc} is solved through the central finite differences
\begin{equation}\label{mesh}
\frac{dx_n}{dt}=\frac{1}{\tau\rho_n \Delta \xi^2}\left( \frac{\rho_{n+1}+\rho_n}{2} (x_{n+1}+x_n) - \frac{\rho_{n}+\rho_{n-1}}{2} (x_{n}+x_{n-1})\right)
\end{equation}
for $n=1,2,\ldots ,N_I-1$, where the spatial nodes $x_n\in\mathcal{T}_h(t)$ are the unknown solutions of the nodes $\xi_n\in\mathcal{T}_h^c$. The relaxation parameter $\tau$ is specified by the user
for adjusting the response time of mesh movement according to the changes of the density function $\rho(x,t)$. The functions $\rho_{n}$, $n=0,1,\ldots ,N_I$, are computed through the discrete form of the monitor function $\rho$. For instance, we have for the optimal monitor function \eqref{density_opt}
\begin{equation}\label{rho_opt}
\rho_n = \left( 1+\frac{1}{\alpha}\left |u_{xx,n}\right |^2\right)^{1/3}, \; \quad n=0,\ldots , N_I,
\end{equation}
where the terms $u_{xx,n}$ are computed by the central difference approximations using the solutions at the mesh nodes $x_n$. At each time step, further, the discrete monitor functions $\rho_n$ are  smoothed by weighted averaging \cite[Section 1.2]{Huang11amm}.

There are different approaches \cite{Huang11} of the adaptive moving mesh method for solving the fully discrete system \eqref{time}, called the physical PDE, and the MMPDE \eqref{mesh}. One common approach is solving both systems simultaneously using a quasi-Lagrange approach. In this approach, the time derivative term requires a special attention since the mesh is assumed to move in a continuous manner by the time progresses. In such  cases, there occur an extra convective term which may cause difficulties. Additionally, the solution of the both systems simultaneously needs the coupling of the systems, as a result the dimension of the system to be solved increases. Another choice, which we use in this paper, is the alternate solution using a rezoning approach. In this case, the physical PDE and the MMPDE are separated from each other, which allows more flexibility as the mesh generation can be coded  separately and embedded in the solution of the physical PDE. In the static rezoning approach, the change in the spatial mesh is derived in a discrete form similar to the solution \cite{Hyman89,Hyman03,Huang11}. In each time step, first the mesh adaptation is handled using the solution on the old mesh, then the solution is obtained by solving the physical PDE on the newly generated mesh.

In the case of dG method, the unknown coefficient vectors have to satisfy $\upsilon^j\approx \upsilon (t_j)$ in the physical PDE \eqref{time} with the discrete mesh density function \eqref{rho_opt} in the MMPDE \eqref{mesh}. In general, it is difficult to use the coefficient vectors $\upsilon^j$  directly to compute the discrete mesh density functions $\rho_n$, unless we use the Lagrange nodal basis. Let us choose the dG basis functions $\phi_i^n$ as the Lagrange nodal basis functions, $i=1,\ldots ,N_k$, $n=1,\ldots ,N_I$. Then, recalling the definitions \eqref{limit} of the limit terms at a node $x_n$, we have  the relations
\begin{equation}\label{coef}
u_h(x_{n}^-,t)=\upsilon_{N_k}^n(t) \; , \quad u_h(x_{n}^+,t)=\upsilon_{1}^{n+1}(t), \quad n=1,\ldots , N_{I}-1.
\end{equation}
The relations \eqref{coef} originate from the fact that dG methods use discontinuous basis functions at the inter-element nodes. As a result, a dG approximation $u_h(x,t)$ has two traces at an inter-element node $x_n$ from the neighboring sub-intervals $I_n$ and $I_{n+1}$, as shown in Figure~\ref{fig_jump}, which are not the same in general. Using the relations \eqref{coef} and recalling again the average definition in \eqref{jump1}, we can accept the value of the approximate solution $u_h(x_n,t)$ as
\begin{equation}\label{val}
\begin{aligned}
u_h(x_{n},t):&= \{ u_h(x_{n},t)\} = \frac{1}{2} (u_h(x_{n}^-,t)+u_h(x_{n}^+,t))\\
&=\frac{1}{2} (\upsilon_{N_k}^n(t)+\upsilon_{1}^{n+1}(t)), \quad n=1,\ldots , N_{I}-1.
\end{aligned}
\end{equation}
In this way, the solutions  of the physical PDE \eqref{time} will be  consistent with the MMPDE \eqref{mesh}, where the computation of the discrete mesh density functions $\rho_n$ require the discrete approximations $u_h(x,t)$ at the inter-element nodes $x_n$, $n=0,\ldots ,N_I$.

\begin{algorithm}
Given the current spatial mesh $\mathcal{T}_h(t_{j-1}): \; x_0^{j-1}<x_1^{j-1}<\ldots < x_{N_I}^{j-1}$, coefficient vector $\upsilon^{j-1}$, the parameter $\tau$ and time-step size $\Delta t$, do for $j=1,\ldots , J$,
\begin{algorithmic}[1]
\STATE  Compute a temporary coefficient vector $\tilde{\upsilon}^j$ by solving the physical PDE \eqref{time} on the current mesh $\mathcal{T}_h(t_{j-1})$,
\STATE  According to \eqref{val}, calculate the consistent values of $\tilde{u}_h(x_n,t_j)$ using $\tilde{\upsilon}^j$,
\STATE  Compute the discrete mesh density functions $\tilde{\rho}_n$ using the values $\tilde{u}_h(x_n,t_j)$,
\STATE  Find the mesh $\mathcal{T}_h(t_{j}): \; x_0^{j}<x_1^{j}<\ldots < x_{N_I}^{j}$ by solving the MMPDE \eqref{mesh} for the discrete mesh density functions $\tilde{\rho}_n$,
\STATE  Interpolate the coefficient vector $\upsilon^{j-1}$ to be used in the new mesh $\mathcal{T}_h(t_{j})$,
\STATE  Compute the coefficient vector $\upsilon^j$ by solving the physical PDE \eqref{time} on the new mesh $\mathcal{T}_h(t_{j})$,
\STATE  Go to next time step.
\end{algorithmic}
\caption{Moving mesh algorithm\label{alg}}
\end{algorithm}

The general algorithm is summarized in Algorithm~\ref{alg}. We start with an initial mesh $\mathcal{T}_h(t_0)$ (possibly a uniform mesh). Then on an arbitrary time-step $(t_{j-1},t_j]$, we, firstly, solve the physical PDE \eqref{time} on the mesh $\mathcal{T}_h(t_{j-1})$ for an auxiliary solution $\tilde{\upsilon}^j$. After, following the relation \eqref{val}, we use the solution $\tilde{\upsilon}^j$ in the MMPDE \eqref{mesh} to obtain the new mesh $\mathcal{T}_h(t_j)$. Finally, we solve the physical PDE \eqref{time} on the mesh $\mathcal{T}_h(t_{j})$ for the solution $\upsilon^j$, and we proceed to the next time-step. On the other hand, the known solution vector $\upsilon^{j-1}$ will be no more consistent with the new mesh $\mathcal{T}_h(t_j)$. This is a natural consequence of the alternate solution with rezoning approach. To make the known solution $\upsilon^{j-1}$ consistent with the updated spatial mesh $\mathcal{T}_h(t_j)$, we interpolate it between the meshes $\mathcal{T}_h(t_{j-1})$ and $\mathcal{T}_h(t_{j})$. The interpolation procedure is as follows: let $\mathcal{T}_h(t_{j-1}): \; x_0^{j-1}<x_1^{j-1}<\ldots < x_{N_I}^{j-1}$ be the current mesh and $\mathcal{T}_h(t_{j}): \; x_0^{j}<x_1^{j}<\ldots < x_{N_I}^{j}$ denote the updated mesh. For the sake of simplicity, let also $\hat{\upsilon} :=\upsilon^{j-1}$ with entries $\hat{\upsilon}_i^n$ being the coefficient vector defined on the current mesh $\mathcal{T}_h(t_{j-1})$, and $\upsilon$ with entries $\upsilon_i^n$ denotes the interpolated coefficient vector of $\hat{\upsilon}$ into the new mesh $\mathcal{T}_h(t_{j})$, $i=1,\ldots , N_k$, $n=1,\ldots ,N_I$. The local nature of dG methods leads for any $x\in I_s=[x_{s-1}^{j-1},x_{s}^{j-1}]$ to
\begin{equation}\label{local}
u_h(x,t_{j-1})= \sum_{i=1}^{N_k} \hat{\upsilon}_i^s \phi_i^s(x), \quad s=1,\ldots ,N_I.
\end{equation}

On the other hand, using the  Lagrange basis functions, on an arbitrary element we obtain $I_n=[x_{n-1}^{j},x_{n}^{j}]$ on the new mesh $\mathcal{T}_h(t_{j})$, uniformly distributed $N_k$ nodal degrees of freedoms $x_{n,d}^{j}\in I_n$ such that
\begin{equation}\label{dof}
x_{n,d}^{j}= x_{n-1}^{j} + (d-1)\tau ,\quad u_h(x_{n,d}^{j},t_{j-1}) \approx \upsilon_d^n,
\end{equation}
for $d=1,\ldots ,N_k$, and with $\tau =(x_{n}^{j}-x_{n}^{j-1})/(N_k-1)$. In other words, the entries $\upsilon_i^n$ of the interpolated coefficient vector $\upsilon$ are the approximate function values $u_h(x,t_{j-1})$ at the nodal degrees of freedoms $x_{n,d}^j$, $n=1,\ldots ,N_I$, $d=1,\ldots ,N_k$.  For any nodal degrees of freedom $x_{n,d}^j$ on the new mesh $\mathcal{T}_h(t_{j})$, we have to determine the intervals $I_s=[x_{s-1}^{j-1},x_s^{j-1}]$ on the current mesh $\mathcal{T}_h(t_{j-1})$ such that $x_{n,d}^j\in I_s$. Then, using the expansion \eqref{local}, we will be able to obtain the entries of $\upsilon$ as $\upsilon_i^n=u_h(x_{n,i}^j,t_{j-1})$.

\section{Numerical results}
\label{sec:numeric}

In this section we present  several numerical examples demonstrating the effectiveness of the adaptive moving mesh dG method. In all of the examples, the traveling wave solutions are computed by the optimal mesh density function \eqref{density_opt}, but the corresponding mesh trajectories are given for the optimal \eqref{density_opt}, arc-length \eqref{density_arc} and curvature \eqref{density_curv} mesh density functions. 

\subsection{Burgers' equation}

The first test example is the Burgers' equation \cite{Huang11,Soheili12anm}
\begin{equation*}
u_t = \epsilon u_{xx} - \frac{\partial }{\partial  x}\left( \frac{1}{2}u^2\right)
\end{equation*}
with homogeneous Dirichlet boundary conditions on the space-time domain $(x,t)\in [0, 1]\times [0, 1]$ with the diffusion constant $\epsilon= 10^{-4}$. The initial condition is taken as $u(x,0) = \sin (2\pi x) + 0.5\sin (\pi x)$, and linear dG basis functions are used. We set the time-step size $\Delta t= 0.005$, and we choose the relaxation parameter as $\tau= 10^{-1}$.

\begin{figure}[htb!]
\centering
\includegraphics[width=0.55\textwidth]{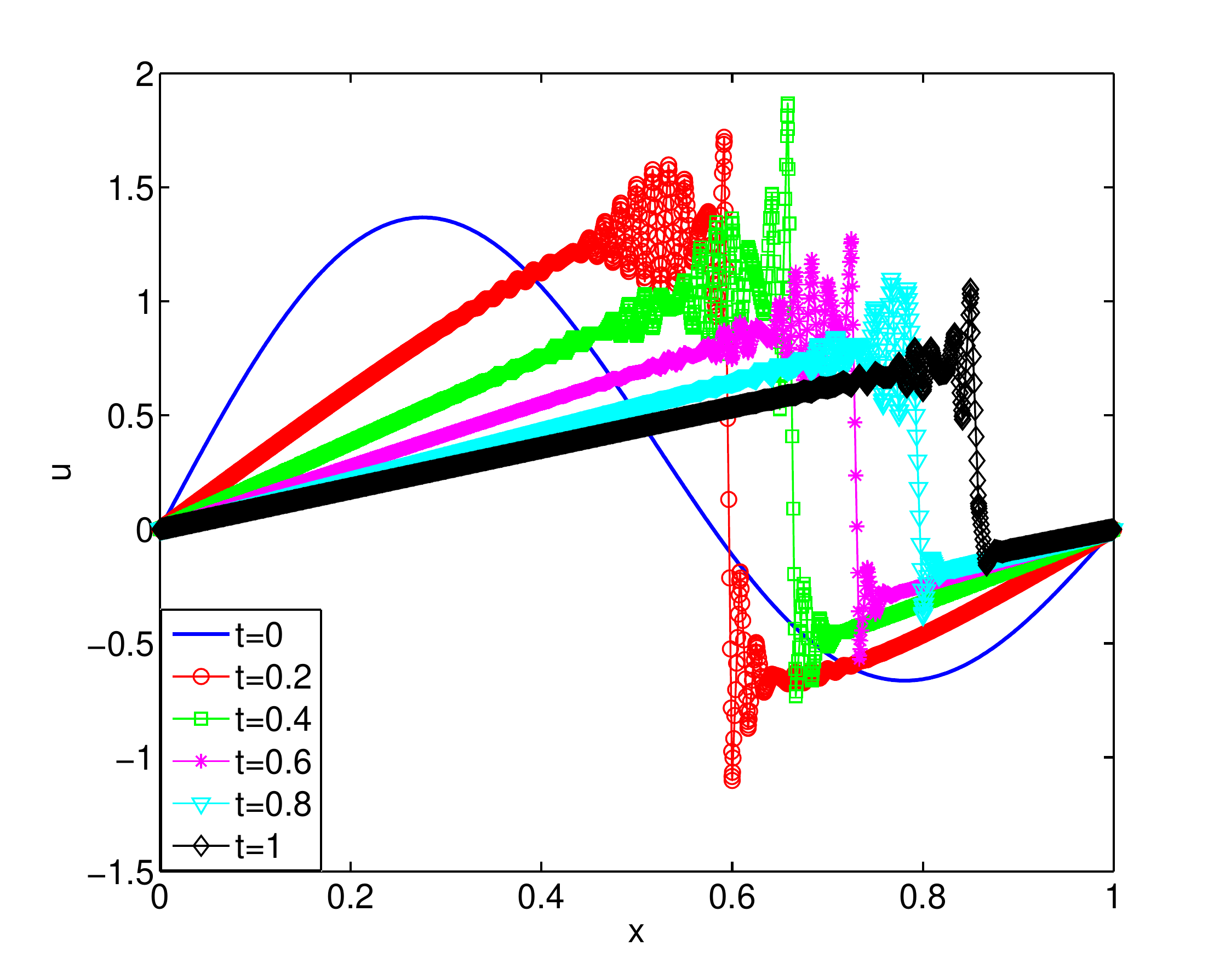}
\caption{Burgers': Solutions at t = 0, 0.2, 0.4, 0.6, 0.8, 1 on the uniform fixed mesh with $N_I=120$ elements.\label{ex1:fix}}
\end{figure}

\begin{figure}[htb]
\centering
\subfloat[(a) Optimal]{\includegraphics[width=0.5\textwidth]{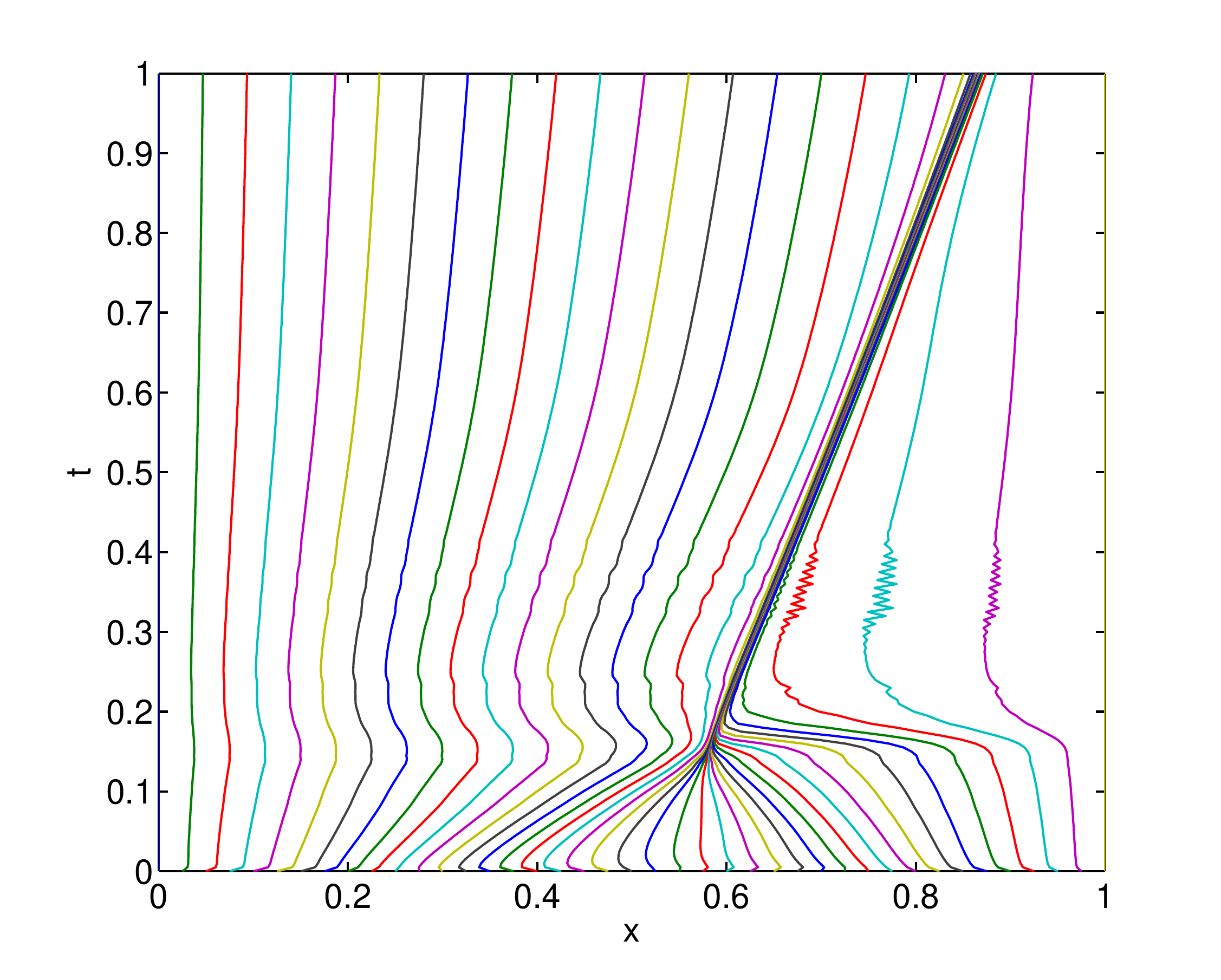}}
\subfloat[(b) Arc-Length]{\includegraphics[width=0.5\textwidth]{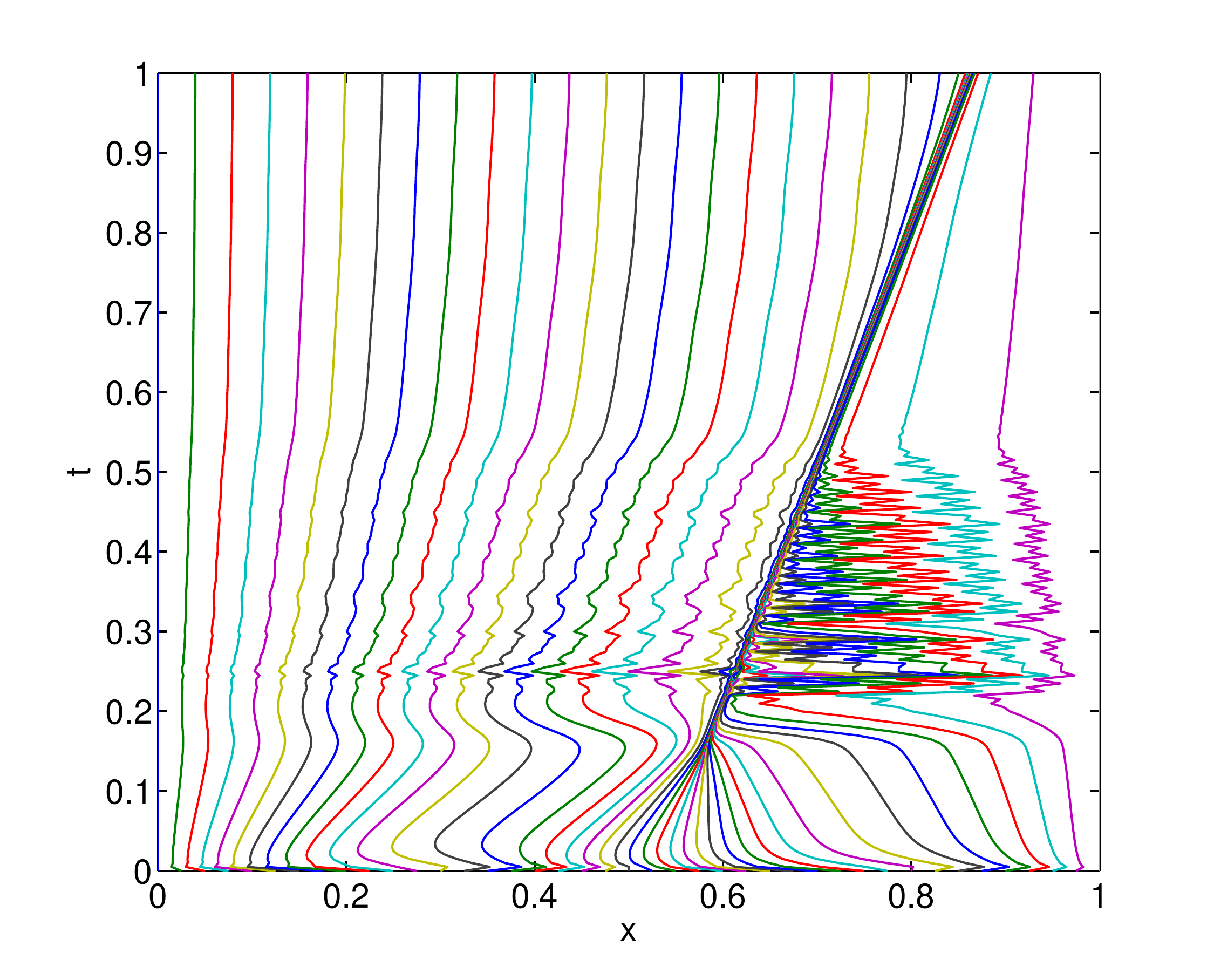}}

\subfloat[(c) Curvature]{\includegraphics[width=0.5\textwidth]{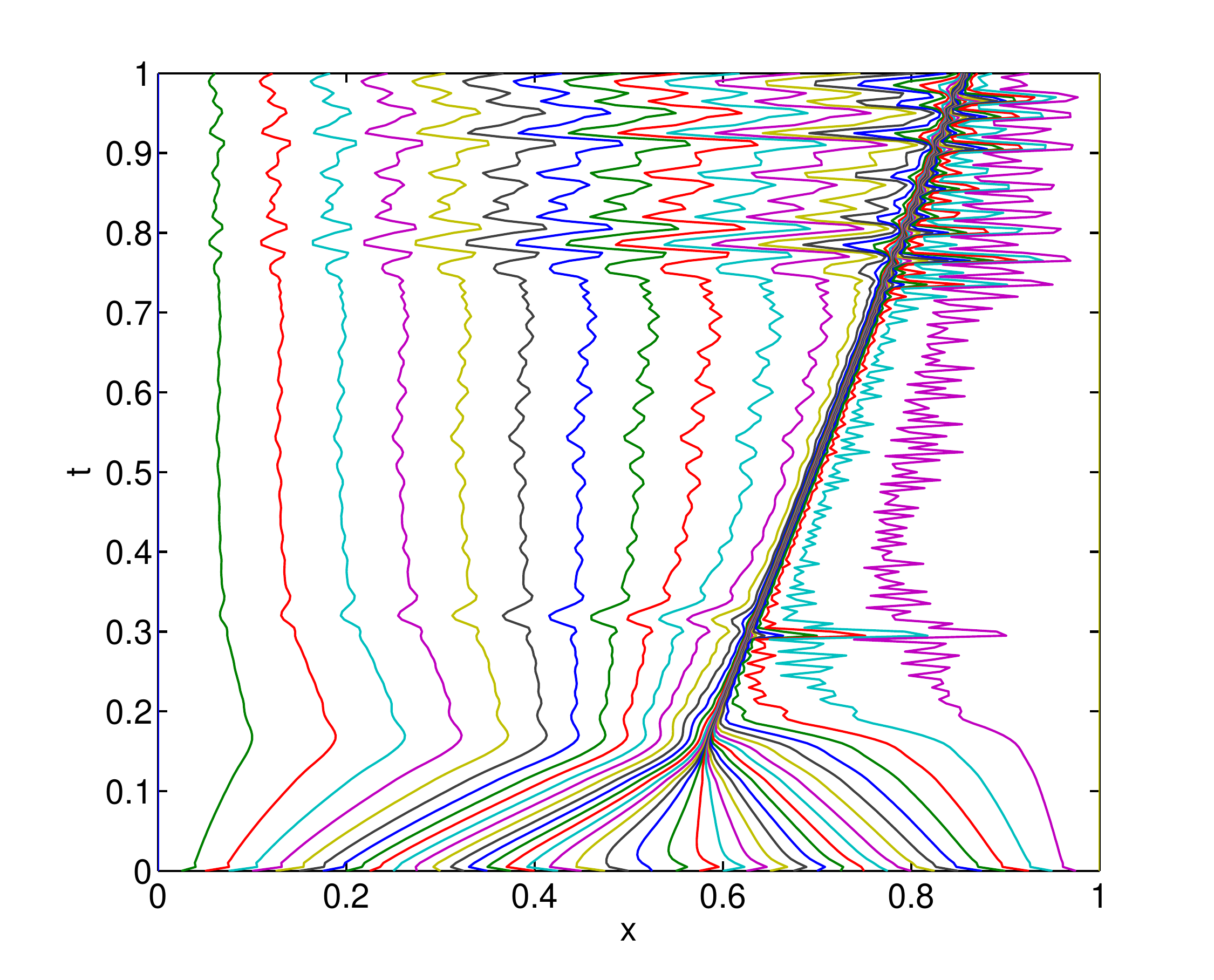}}
\subfloat[(d) Solutions]{\includegraphics[width=0.5\textwidth]{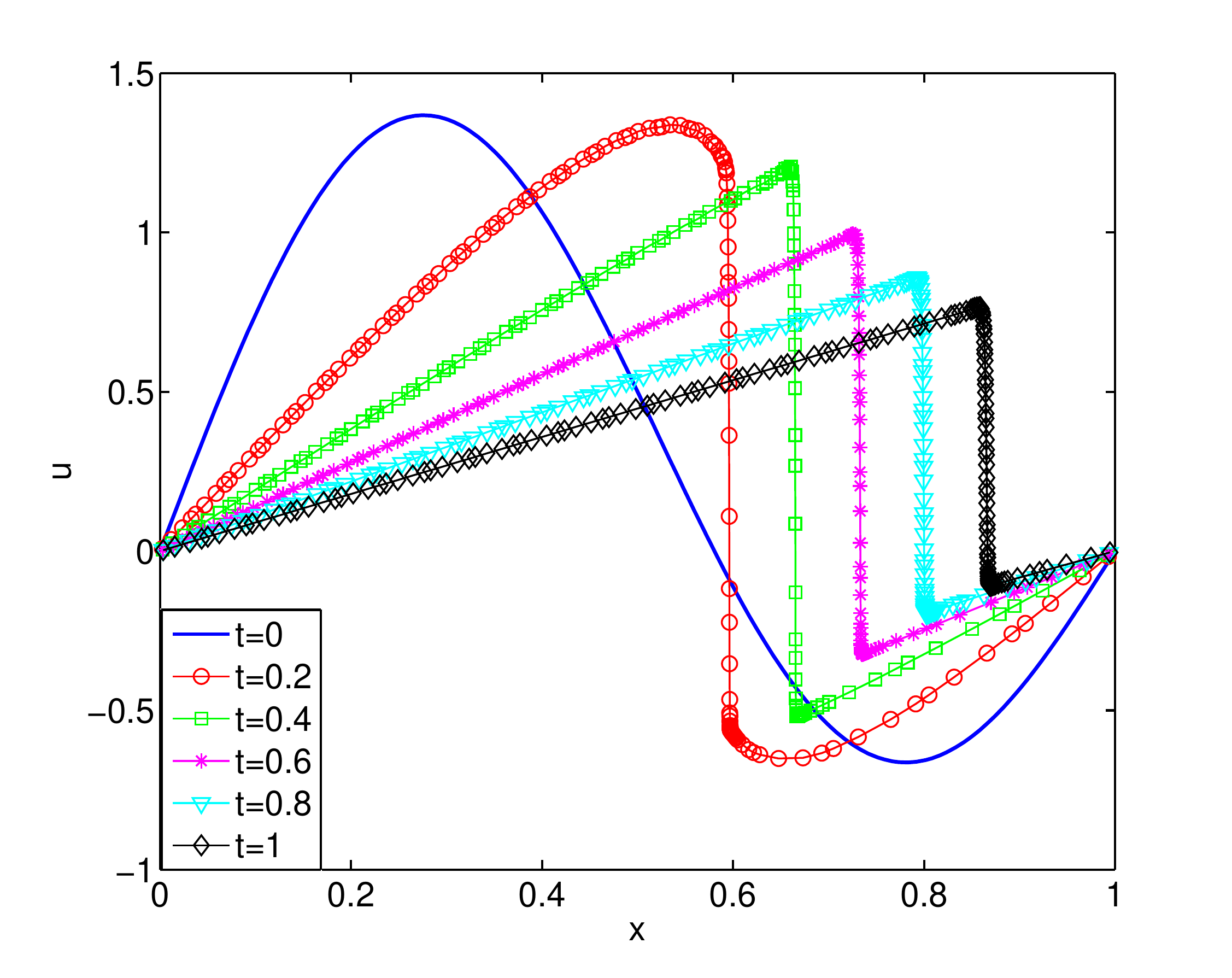}}
\caption{Burgers' equation: (a-c) Moving mesh trajectories with $N_I=40$ elements, (d) solutions at t = 0, 0.2, 0.4, 0.6, 0.8, 1 for the optimal monitor function.\label{ex1:mov}}
\end{figure}

Moving mesh solutions in Figure~\ref{ex1:mov} are capable of  resolving the sharp gradients of the moving fronts in contrast to the oscillatory solutions on the fixed mesh in Figure~\ref{ex1:fix}. In addition, in Figure~\ref{ex1:mov}, mesh trajectories for the curvature and arc-length monitor functions show large variations with respect to time, whereas for the optimal monitor function the mesh trajectories are smooth.

\subsection{Burgers'-Fisher equation}

Consider the Burgers'-Fisher equation \cite{Soheili12anm}
\begin{equation*}
u_t =  u_{xx} -\alpha \frac{\partial }{\partial  x}\left( \frac{1}{2}u^2\right) +  \beta u(u-1),
\end{equation*}
with the exact solution
\begin{equation*}
u(x,t)= \frac{1}{2}\left( 1-\text{tanh}\left( \frac{\alpha}{4}(x-ct)\right)\right),
\end{equation*}
in the space-time domain $(x,t)\in [-1,0]\times (-0.2,0]$. The parameters are $\alpha =24$,  $c= 8$ and $\beta =(2\alpha c-\alpha^2)/4$. We use quadratic dG basis functions and $N_I=40$ spatial elements. The time step-size is taken as $\Delta t=0.001$.

\begin{figure}[htb!]
\centering
\subfloat[(a) Optimal]{\includegraphics[width=0.5\textwidth]{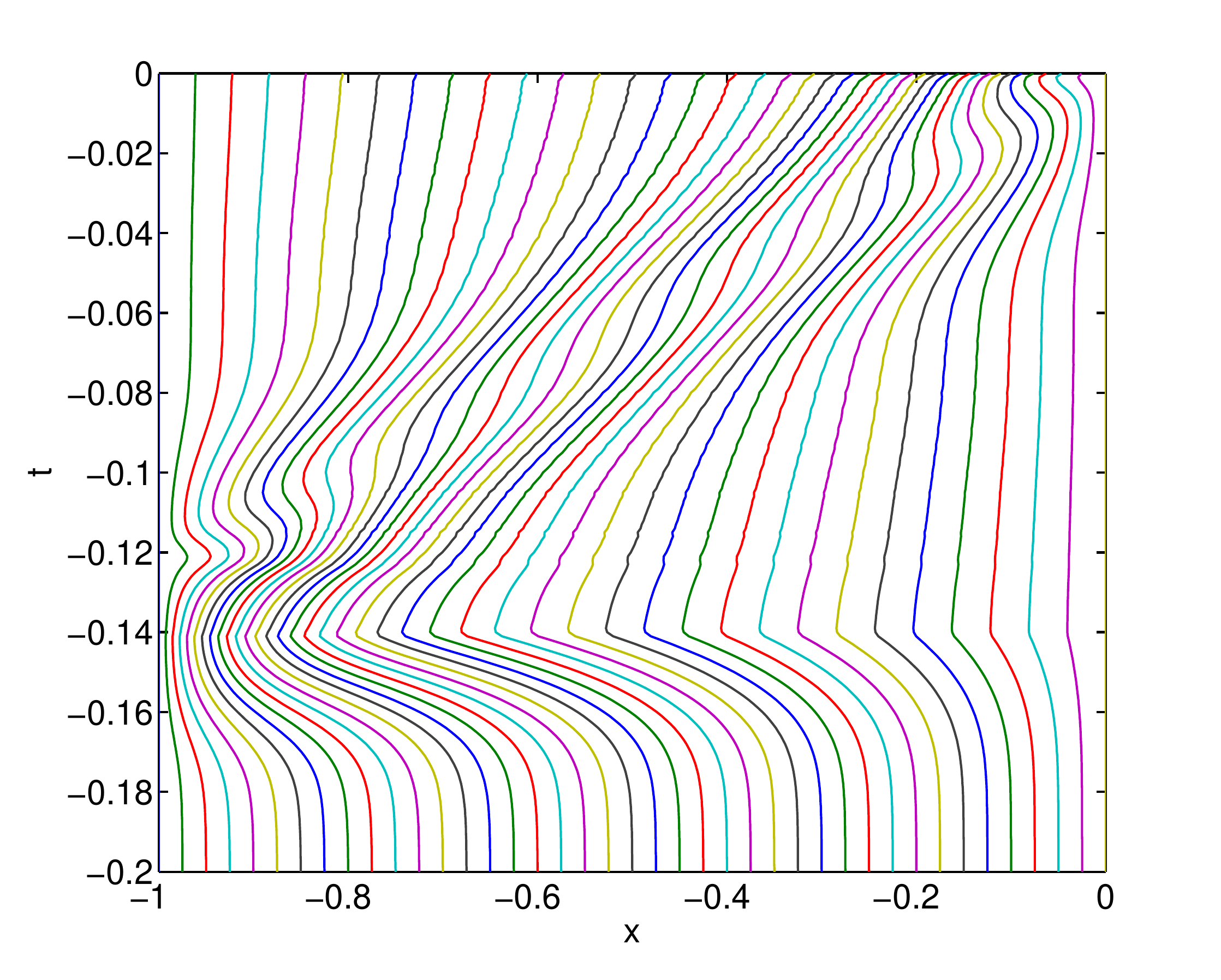}}
\subfloat[(b) Arc-Length]{\includegraphics[width=0.5\textwidth]{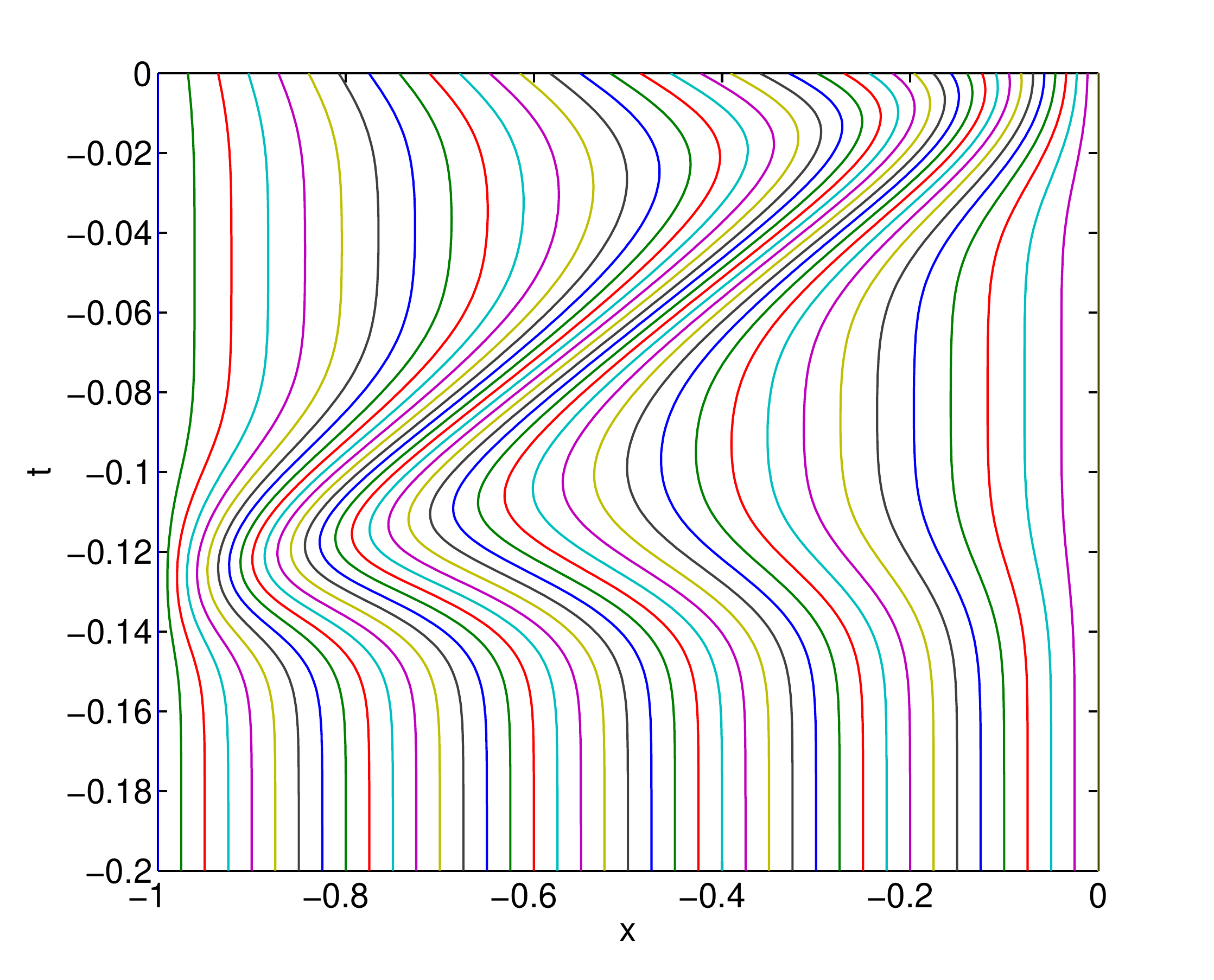}}

\subfloat[(c) Curvature]{\includegraphics[width=0.5\textwidth]{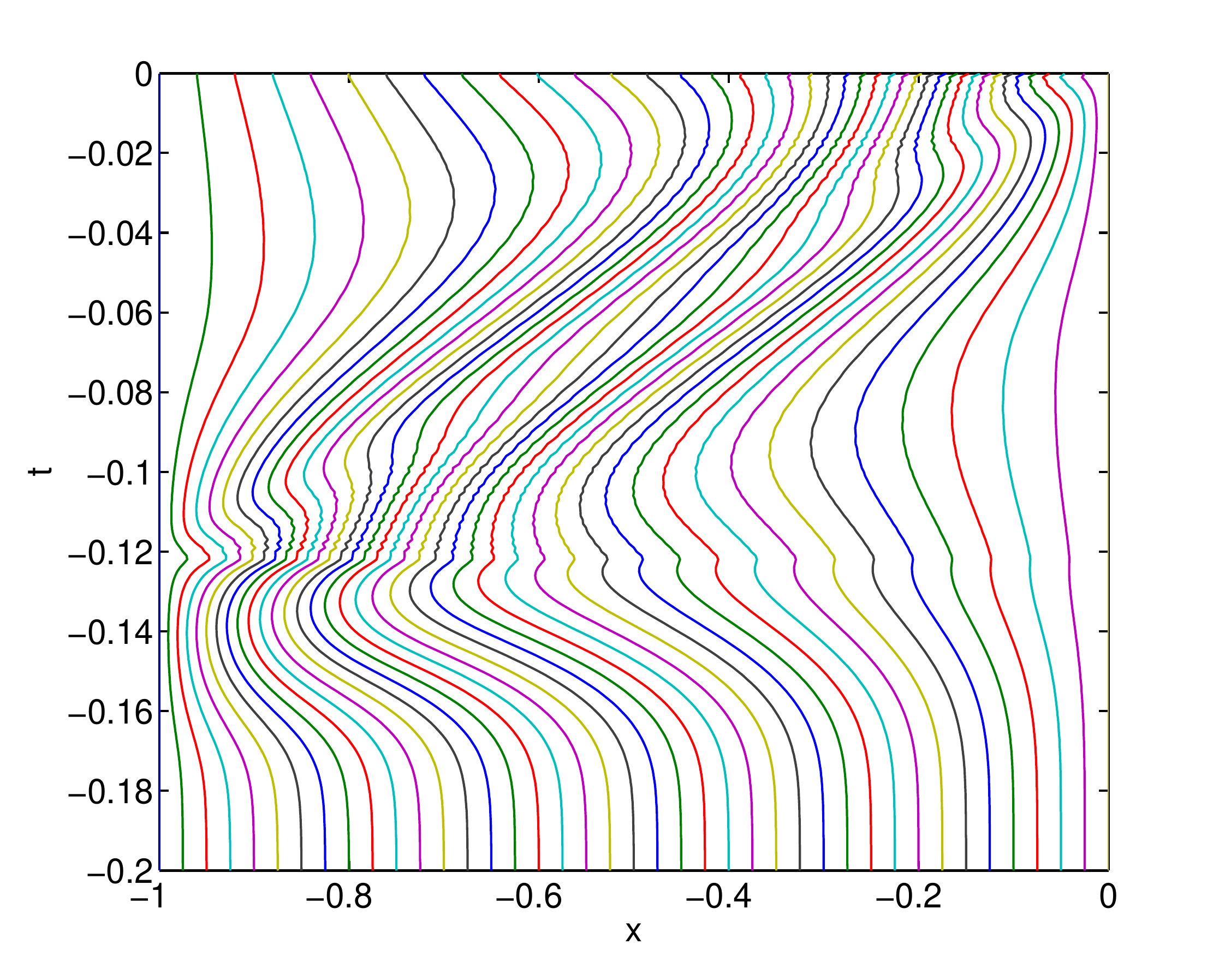}}
\subfloat[(d) Solutions]{\includegraphics[width=0.5\textwidth]{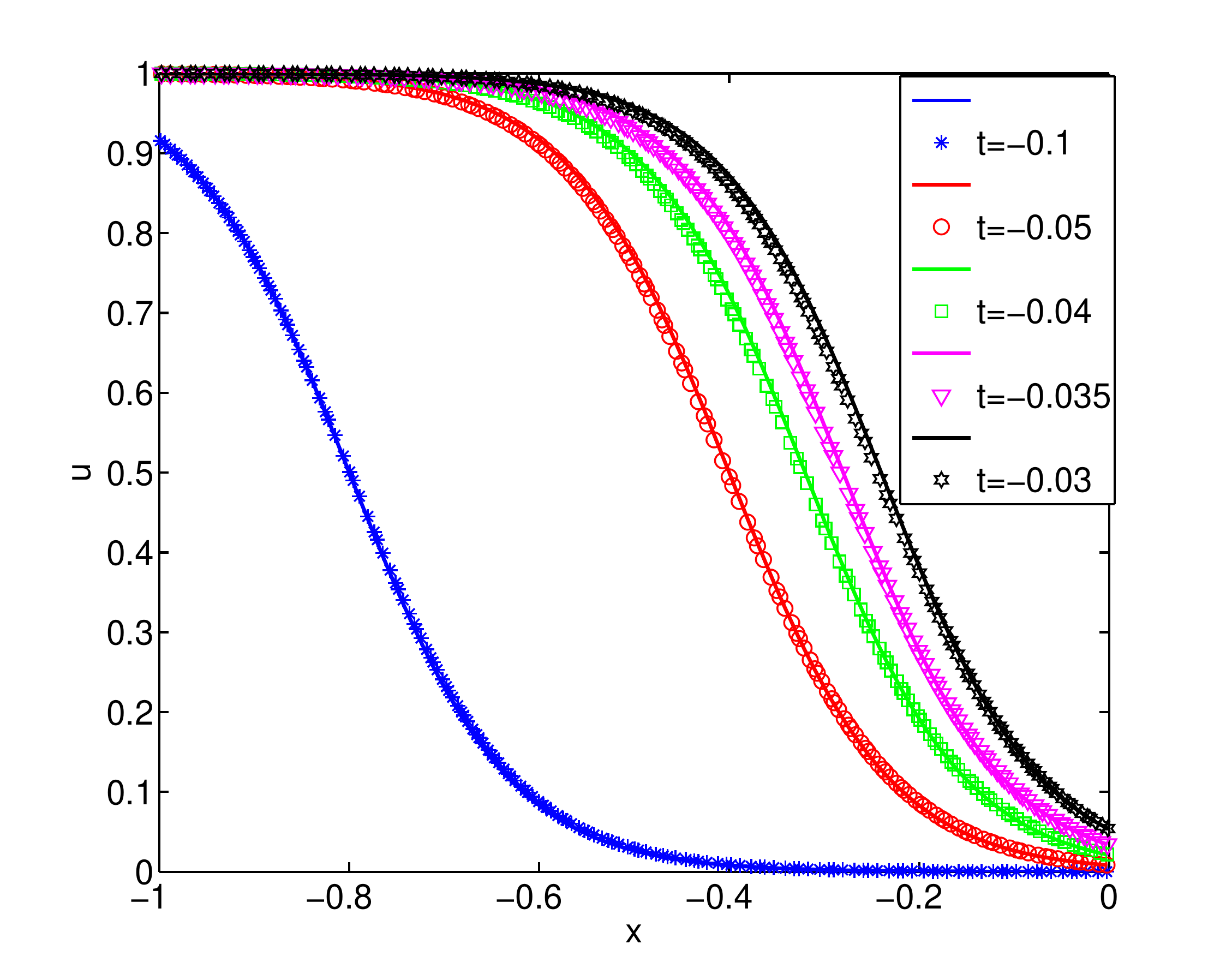}}
\caption{Burgers'-Fisher equation: (a-c) Moving mesh trajectories with $N_I=40$ elements, (d)  solutions at t = -0.1, -0.05, -0.04, -0.035, -0.03 and corresponding exact solutions by solid lines.\label{ex2}}
\end{figure}

In Figure~\ref{ex2}, we give the solutions of Burgers'-Fisher equation at different time instances for the optimal mesh density function, and the moving mesh trajectories obtained by different monitor functions. In Table~\ref{ex2:table} the $L^2$-errors between the exact and numerical solutions are tabulated at different times. The results are very close to those in \cite{Soheili12anm} computed with the same settings.

\begin{table}[htb!]
\centering
\caption{Burgers'-Fisher equation: $L^2$-errors}
\label{ex2:table}
\begin{tabular}{ l | l | l l l l l }
Monitor & $N_I$ & t=-0.1 & t=-0.05  & t=-0.04  & t=-0.035  & t=-0.03 \\
\hline
Optimal  &  40  &  4.6e-03  & 8.4e-03 &  1.1e-02  & 1.1e-02  & 1.3e-02 \\
Arc-Length  &  40  &  2.4e-03  & 4.2e-03 &  5.2e-03  & 5.8e-03 &  6.5e-03 \\
Curvature  &  40  &  2.4e-03  & 4.1e-03 &  5.1e-03 &  5.7e-03 &  6.4e-03 \\
\hline
\end{tabular}
\end{table}

\subsection{Schl\"{o}gl equation}

The final example is Schl\"{o}gl (Nagumo) equation \cite{Buchholz13,Huang11}
\begin{equation*}
u_t =  \epsilon u_{xx} - f(u),
\end{equation*}
with the Ginzburg-Landau free energy
\begin{equation}\label{energy}
\mathcal{E}(u)=\int_{\Omega} \left( \frac{\epsilon}{2}|\nabla u|^2 + F(u)  \right)dx
\end{equation}
with quartic potential function $F(u)=\frac{1}{12\delta} u^2 (3u^2-4(1+\beta)u+6\beta)$ and cubic bi-stable nonlinearity $f(u)=\frac{1}{\delta }u(u-1)(u-\beta)$ satisfying  $f(u)=F^{'}(u)$. 

We consider Schl\"{o}gl equation in the space-time  domain $(x,t)\in (0, 1]\times [0,1]$ with the constant parameter values $\epsilon = 10^{-3}$, $\delta = 10^{-3}$  and $\beta = 0$. The Dirichlet boundary conditions and the initial condition are taken according to the exact solution
\begin{equation*}
u(x,t) = \frac{1}{2}\left( 1-\tanh \left( \frac{x-ct }{\sqrt{8\epsilon \delta}}\right)\right),
\end{equation*}
where $c=\sqrt{\epsilon /2\delta}$ is the speed of the traveling wave.  The solutions are computed with linear and quadratic dG basis functions, and with the time-step size $\Delta t =0.001$. As the number of spatial elements, we take $N_I=120$ and $N_I=40$ to construct a uniform and moving mesh, respectively.

\begin{figure}[htb!]
\centering
\subfloat[(a) Optimal]{\includegraphics[width=0.35\textwidth]{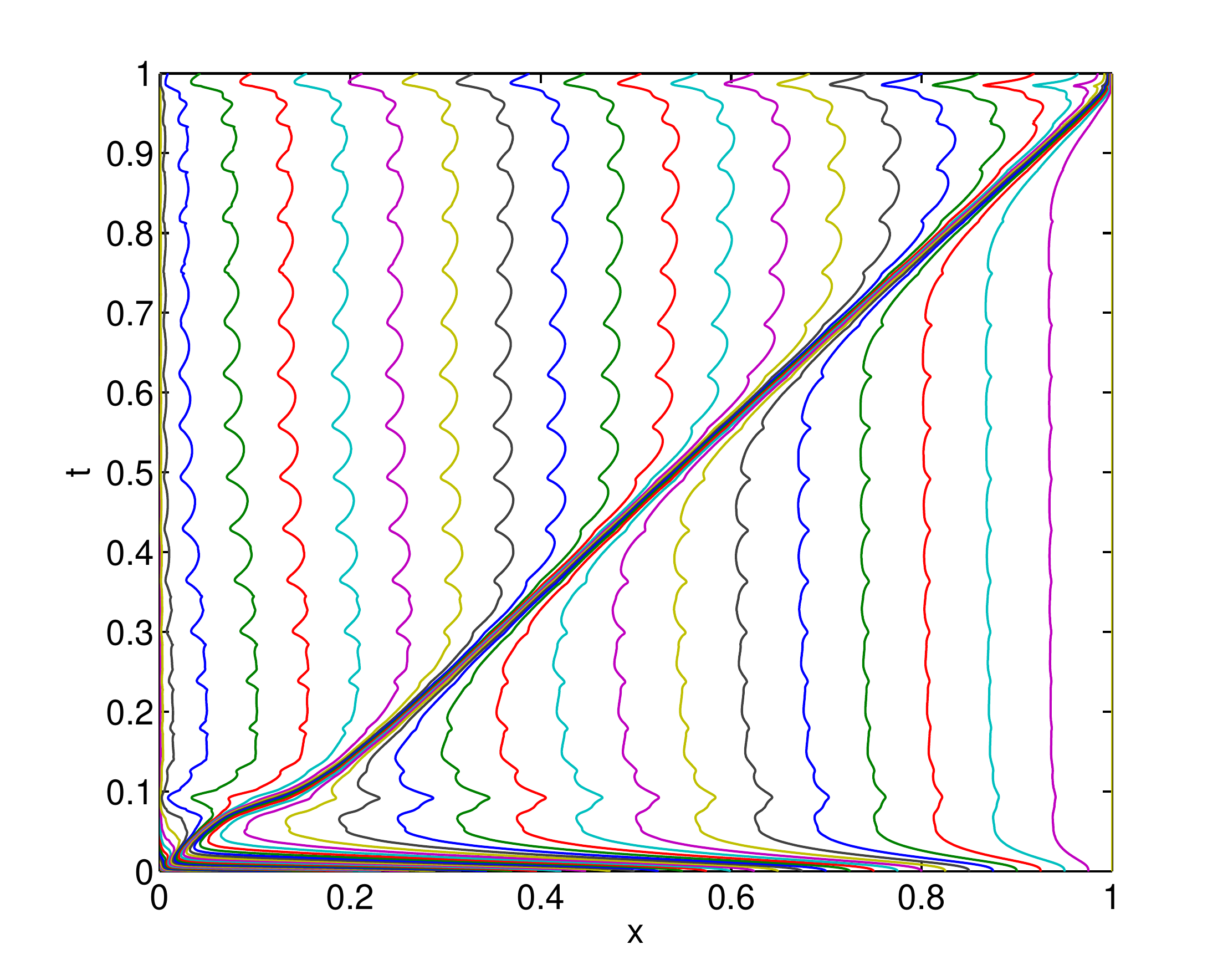}}
\subfloat[(b) Arc-length]{\includegraphics[width=0.35\textwidth]{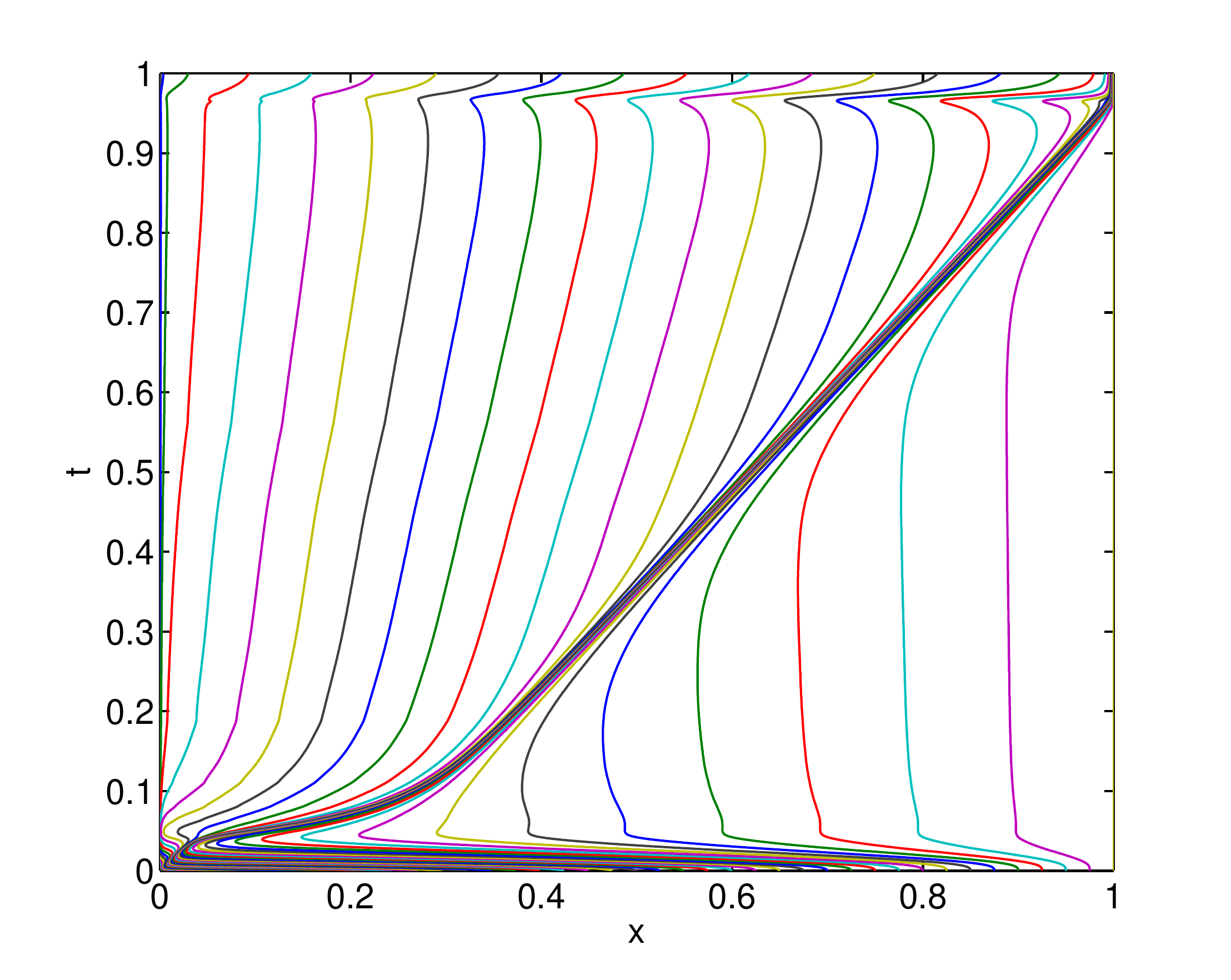}}
\subfloat[(c) Curvature]{\includegraphics[width=0.35\textwidth]{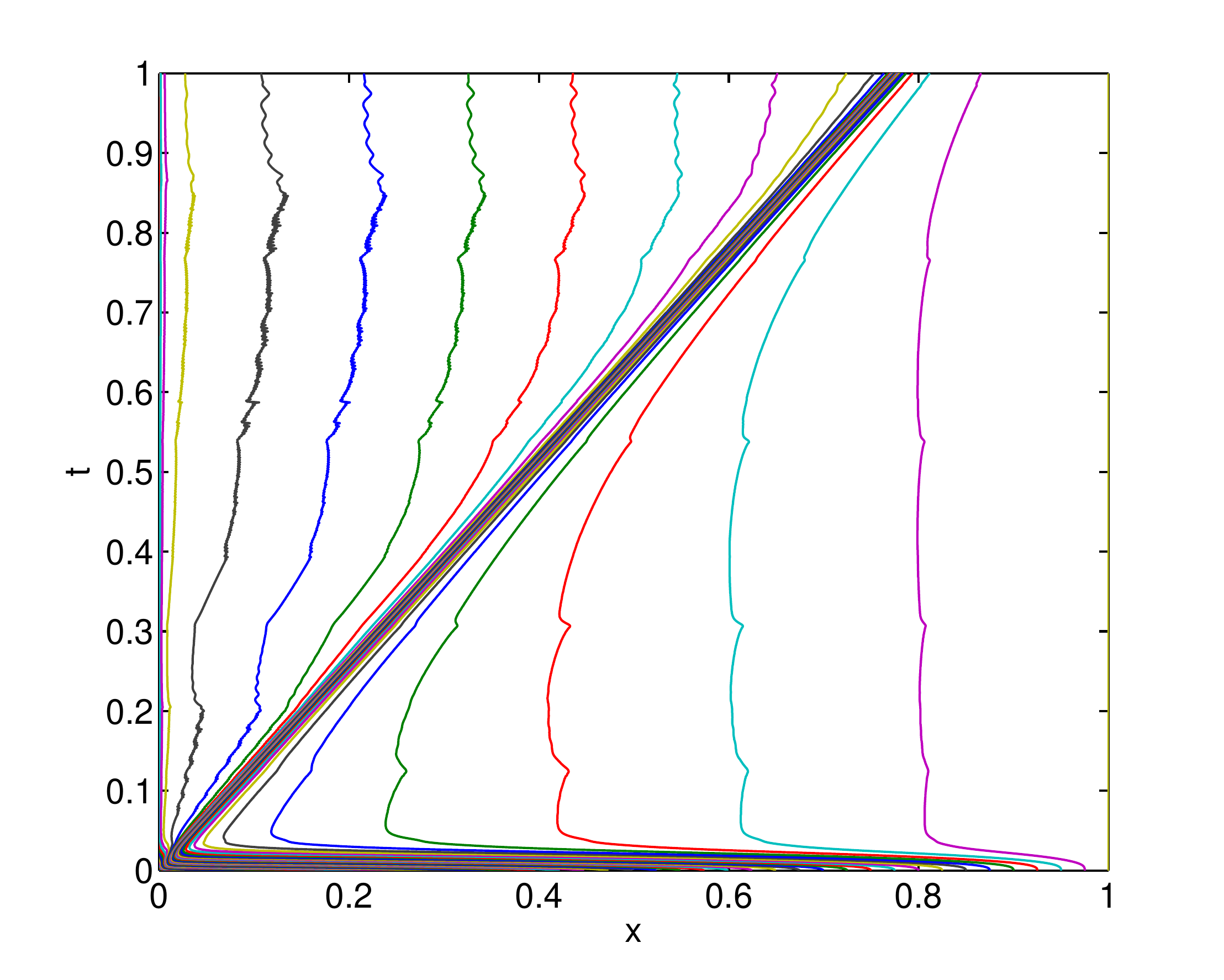}}

\subfloat[(d) Optimal]{\includegraphics[width=0.35\textwidth]{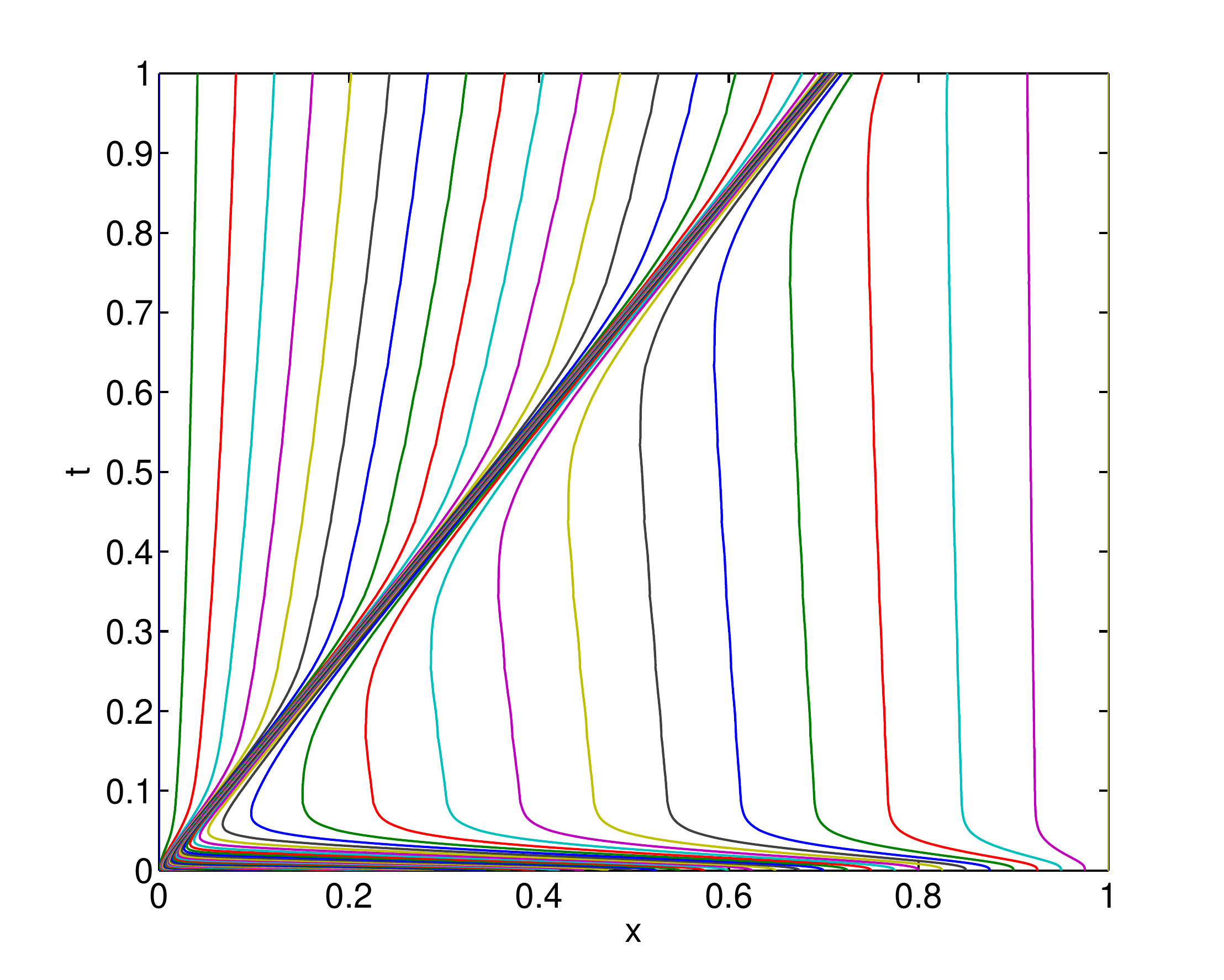}}
\subfloat[(e) Arc-length]{\includegraphics[width=0.35\textwidth]{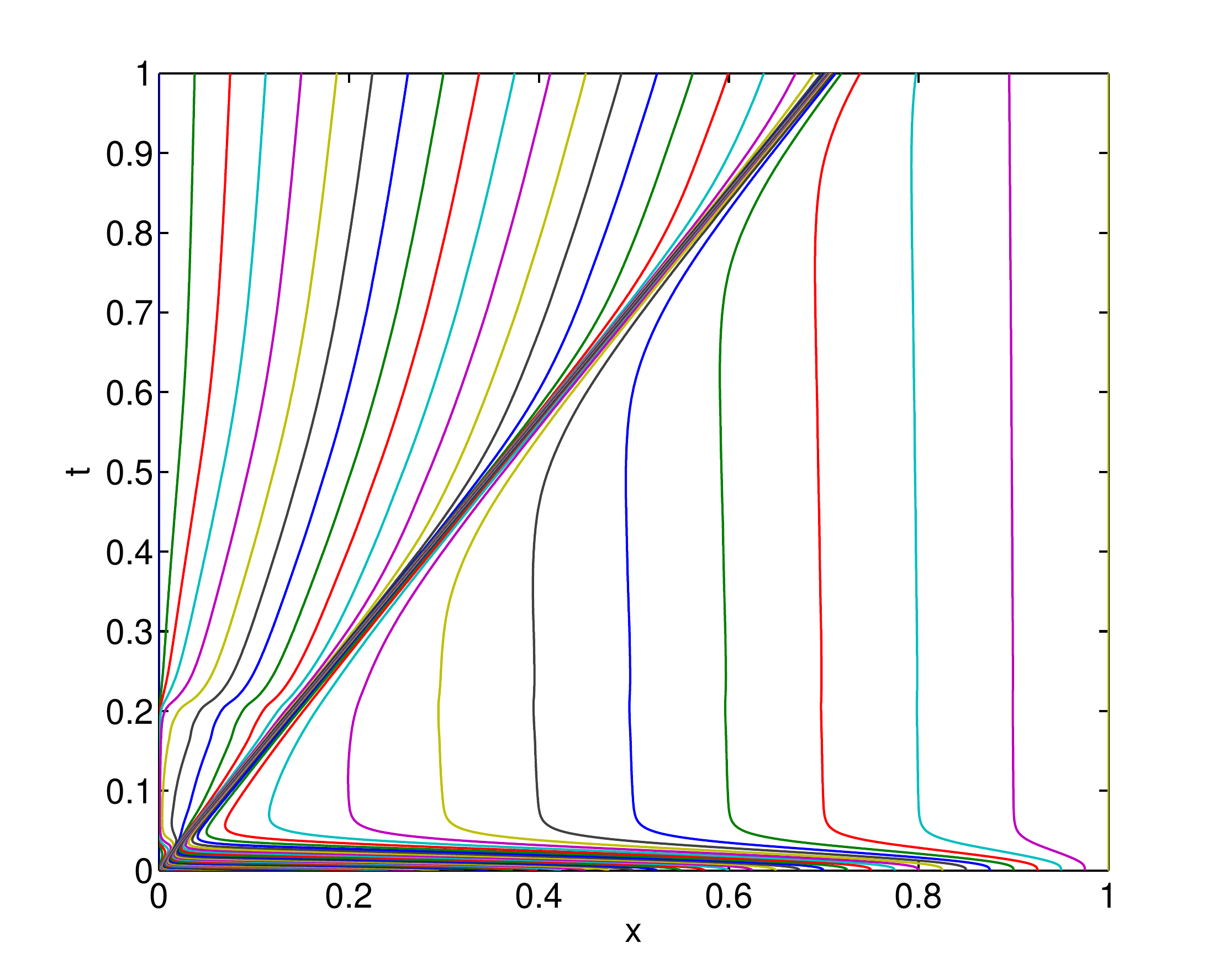}}
\subfloat[(f) Curvature]{\includegraphics[width=0.35\textwidth]{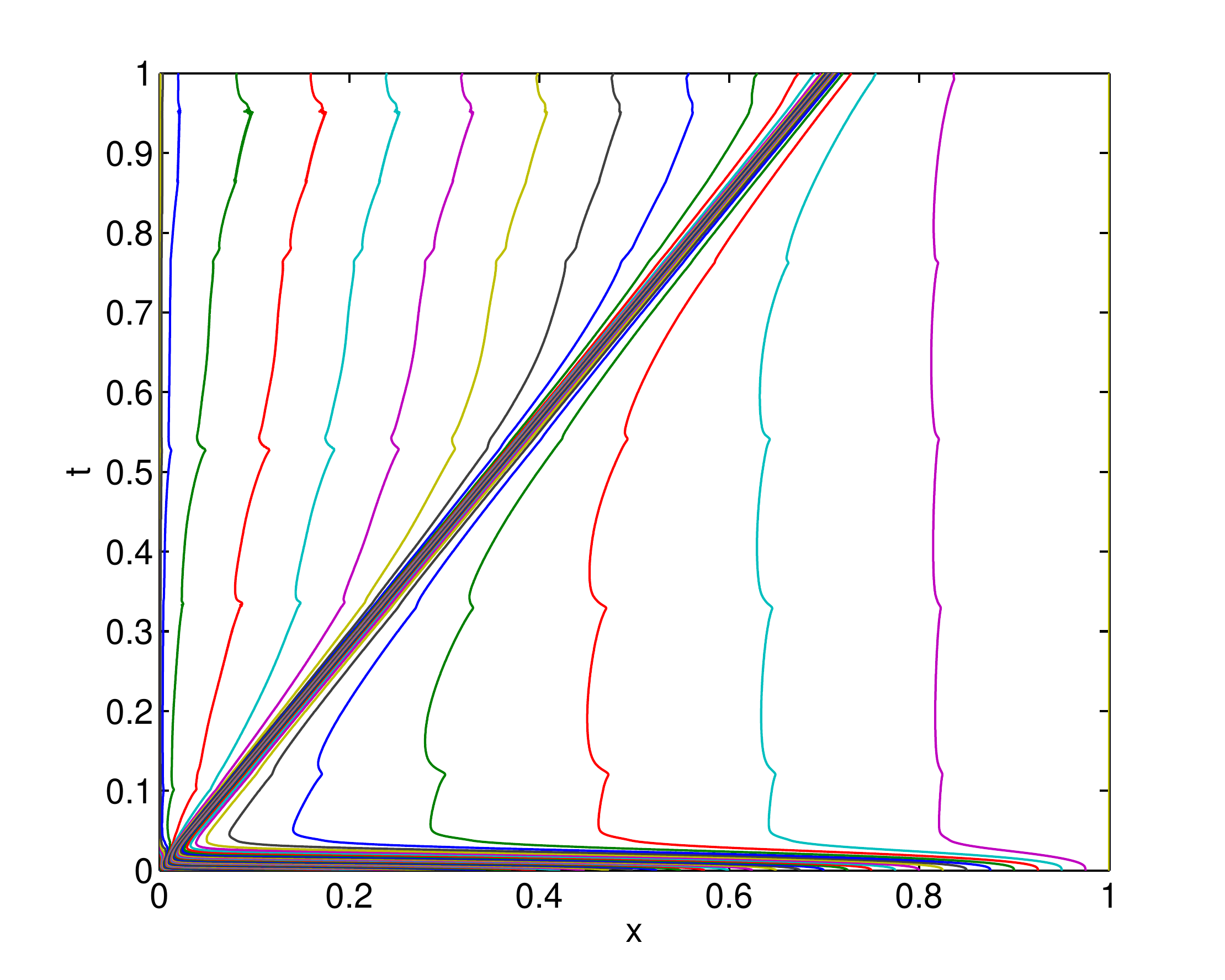}}
\caption{Schl\"{o}gl equation: Moving mesh trajectories with $N_I=40$ elements; (a)-(c) linear dG basis functions, (d)-(f) quadratic dG basis functions.\label{ex3:mesh}}
\end{figure}

\begin{figure}[htb!]
\centering
\subfloat[(a)]{\includegraphics[width=0.37\textwidth]{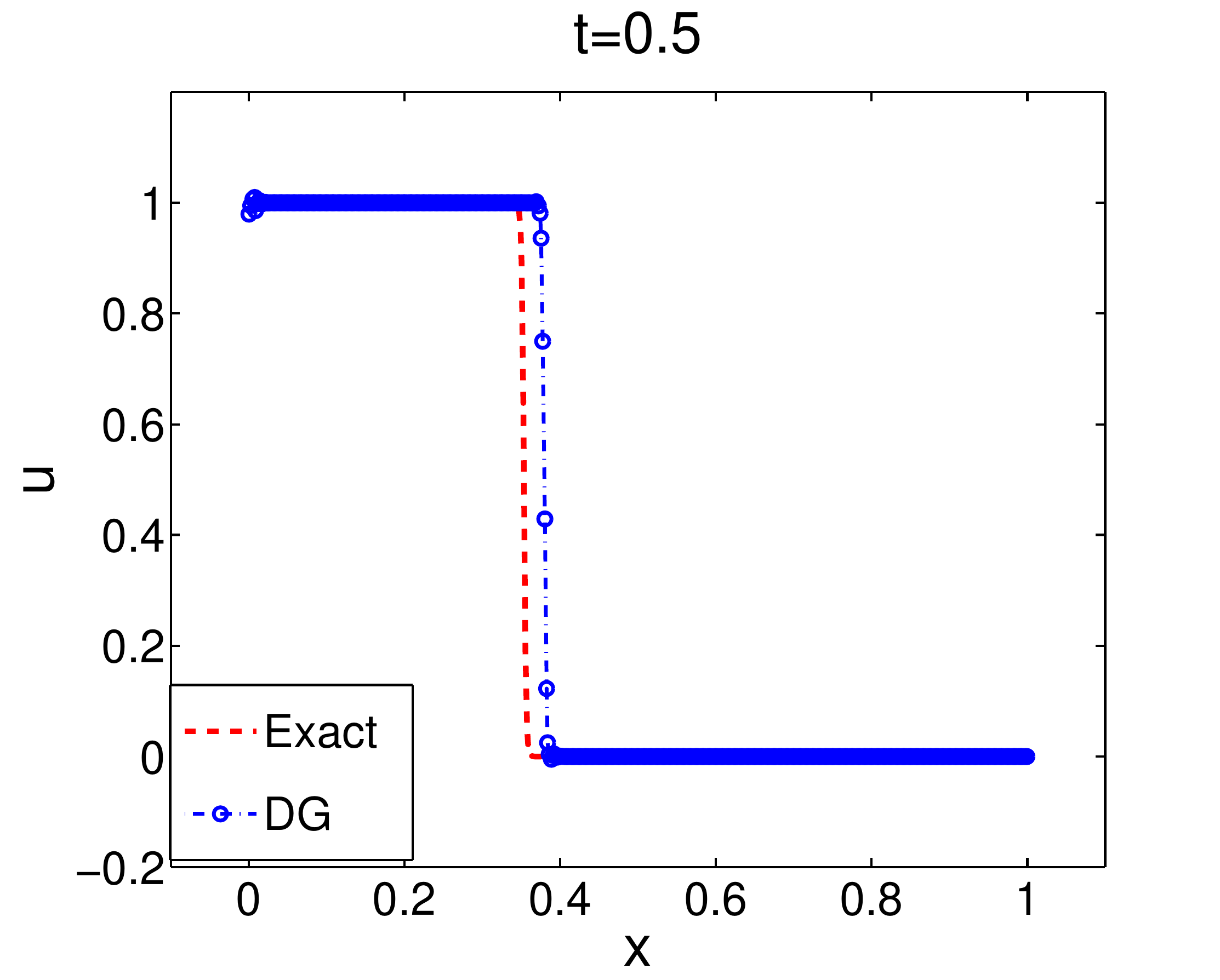}}
\subfloat[(d)]{\includegraphics[width=0.37\textwidth]{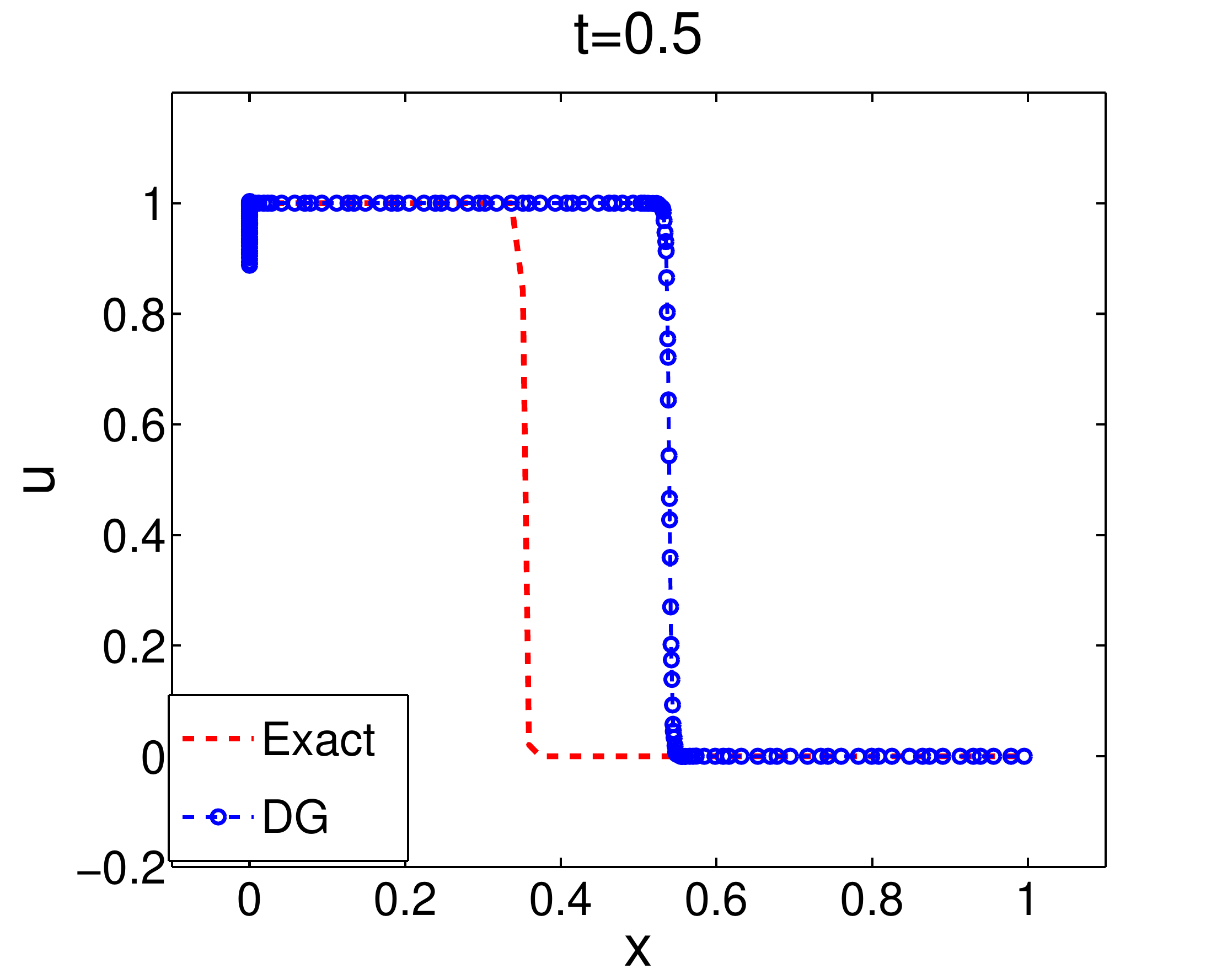}}
\subfloat[(g)]{\includegraphics[width=0.37\textwidth]{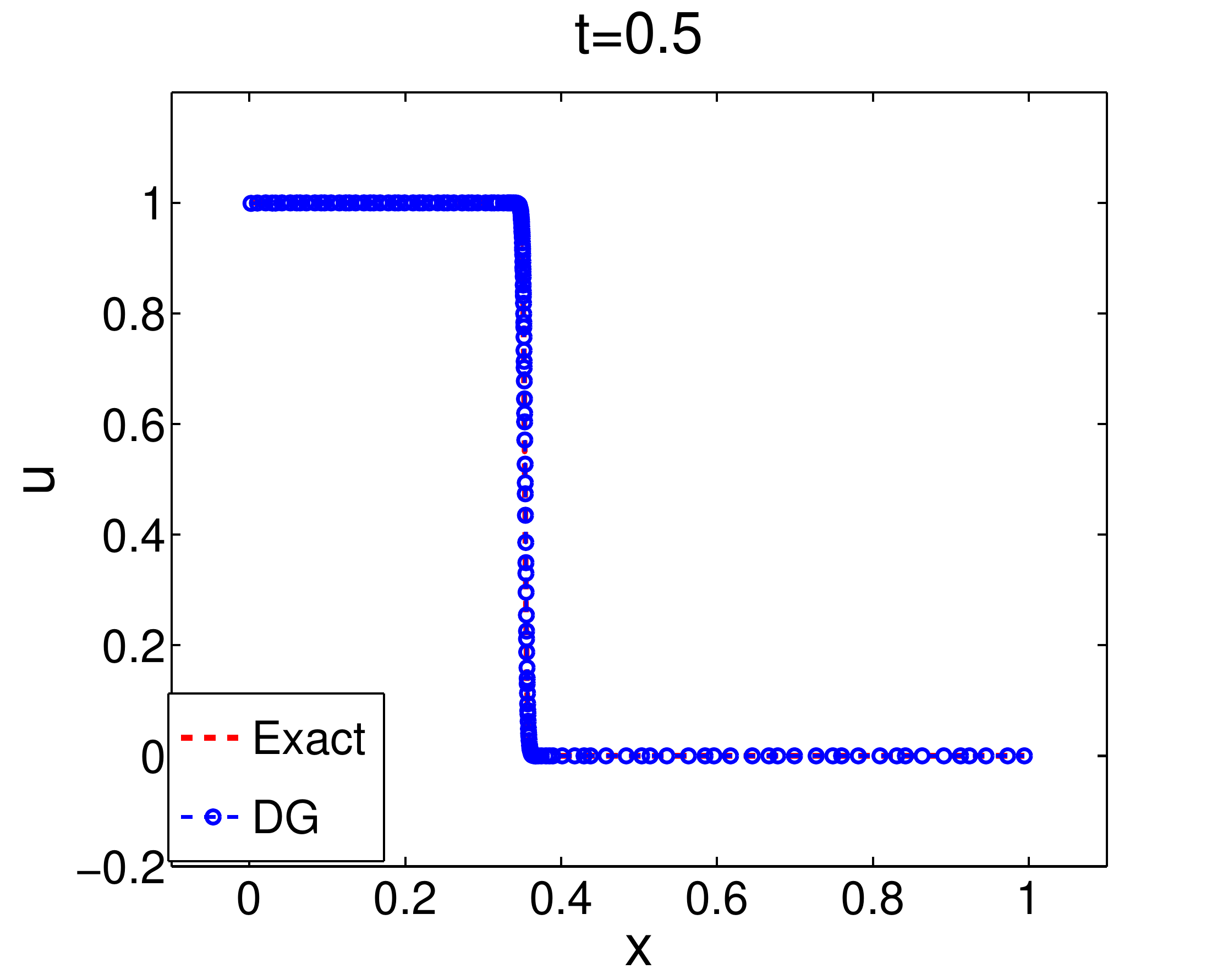}}

\subfloat[(b)]{\includegraphics[width=0.37\textwidth]{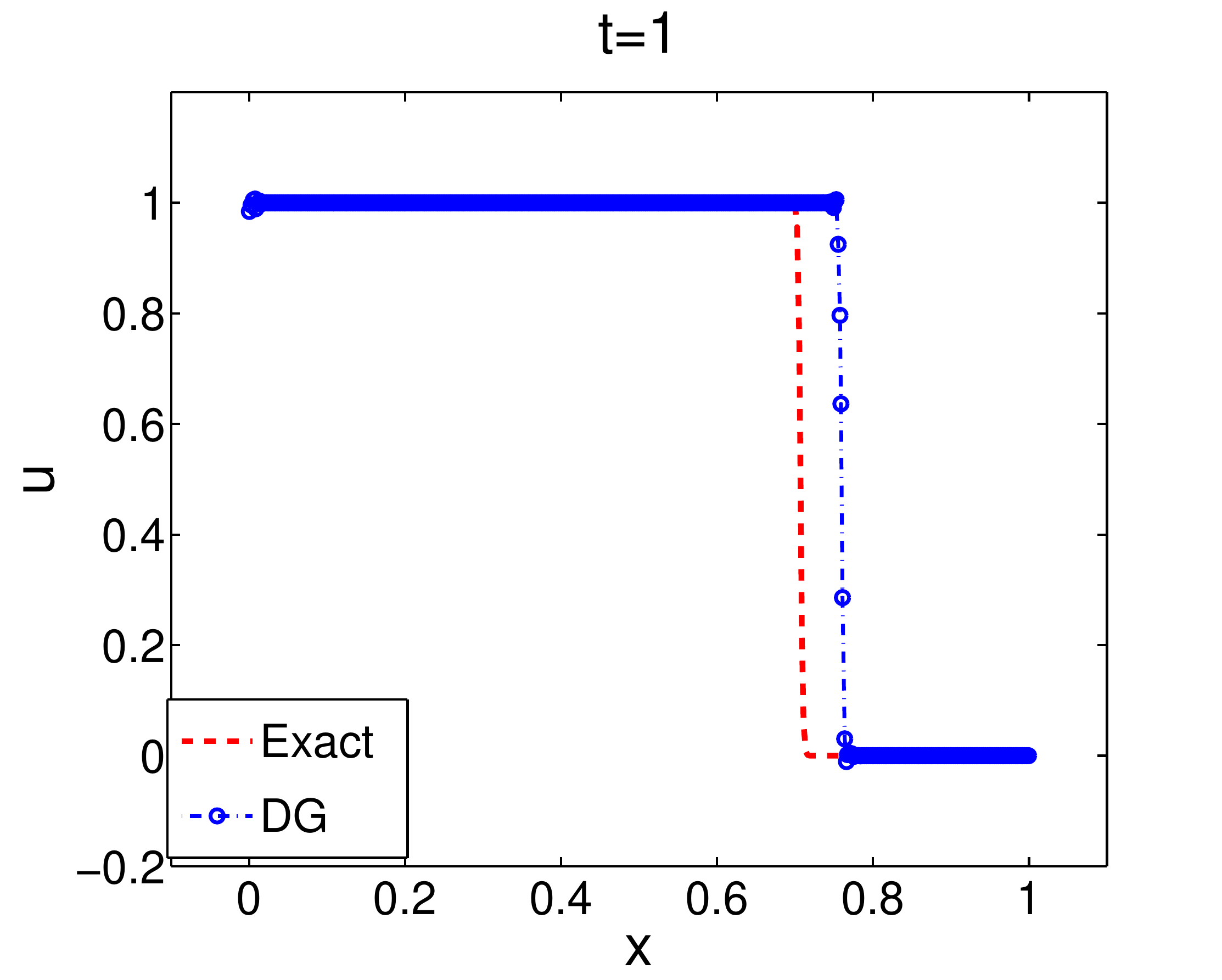}}
\subfloat[(e)]{\includegraphics[width=0.37\textwidth]{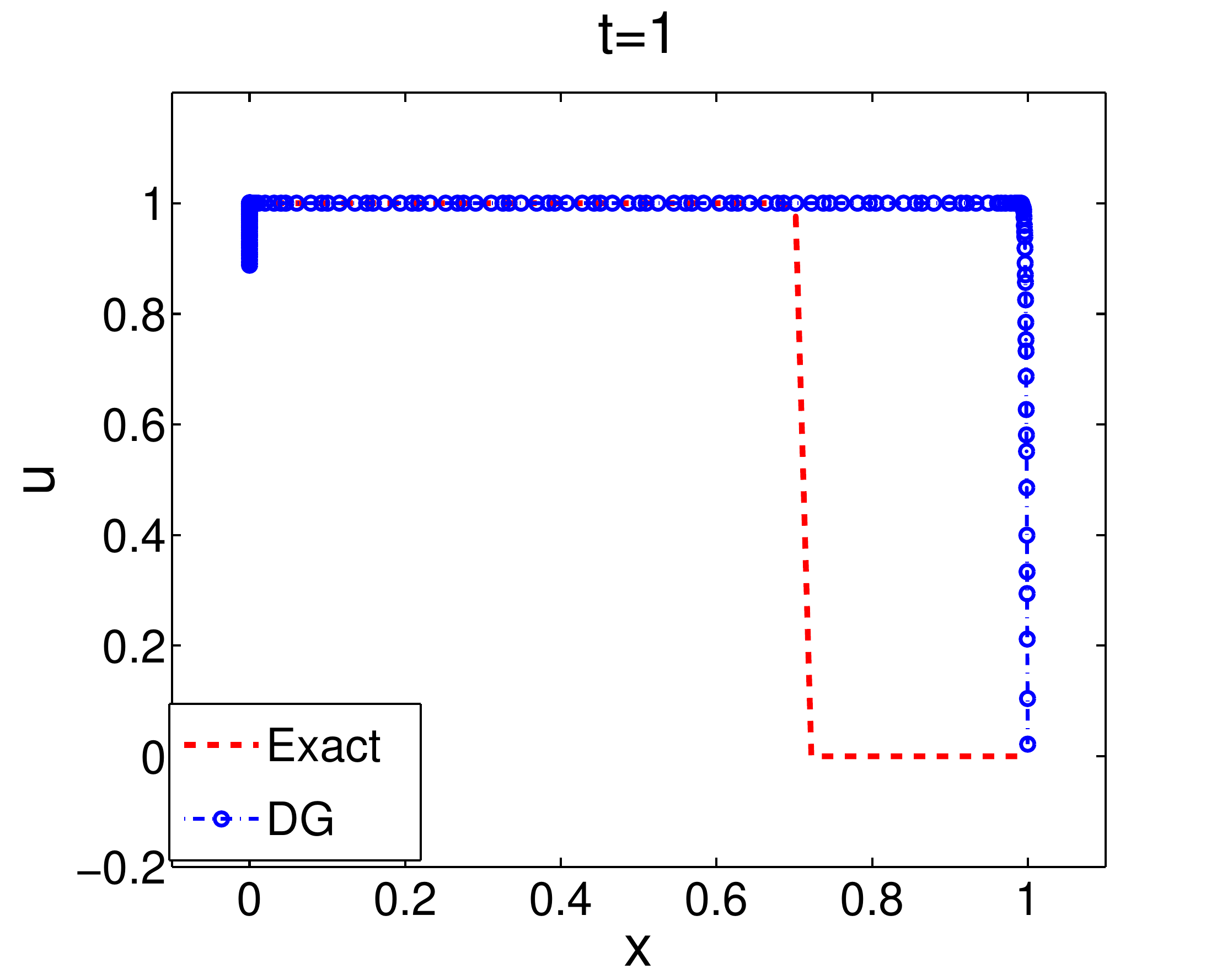}}
\subfloat[(h)]{\includegraphics[width=0.37\textwidth]{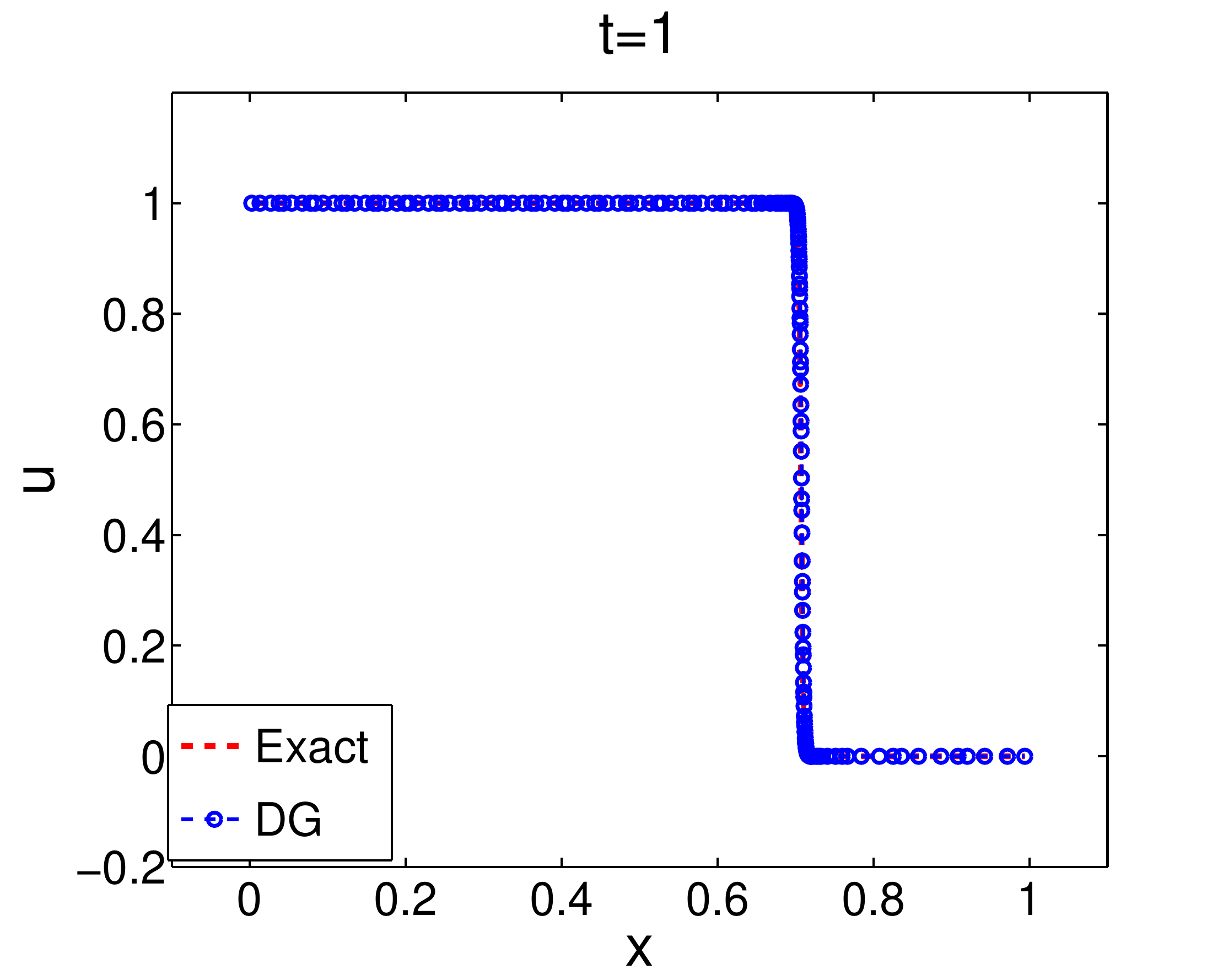}}

\subfloat[(c)]{\includegraphics[width=0.37\textwidth]{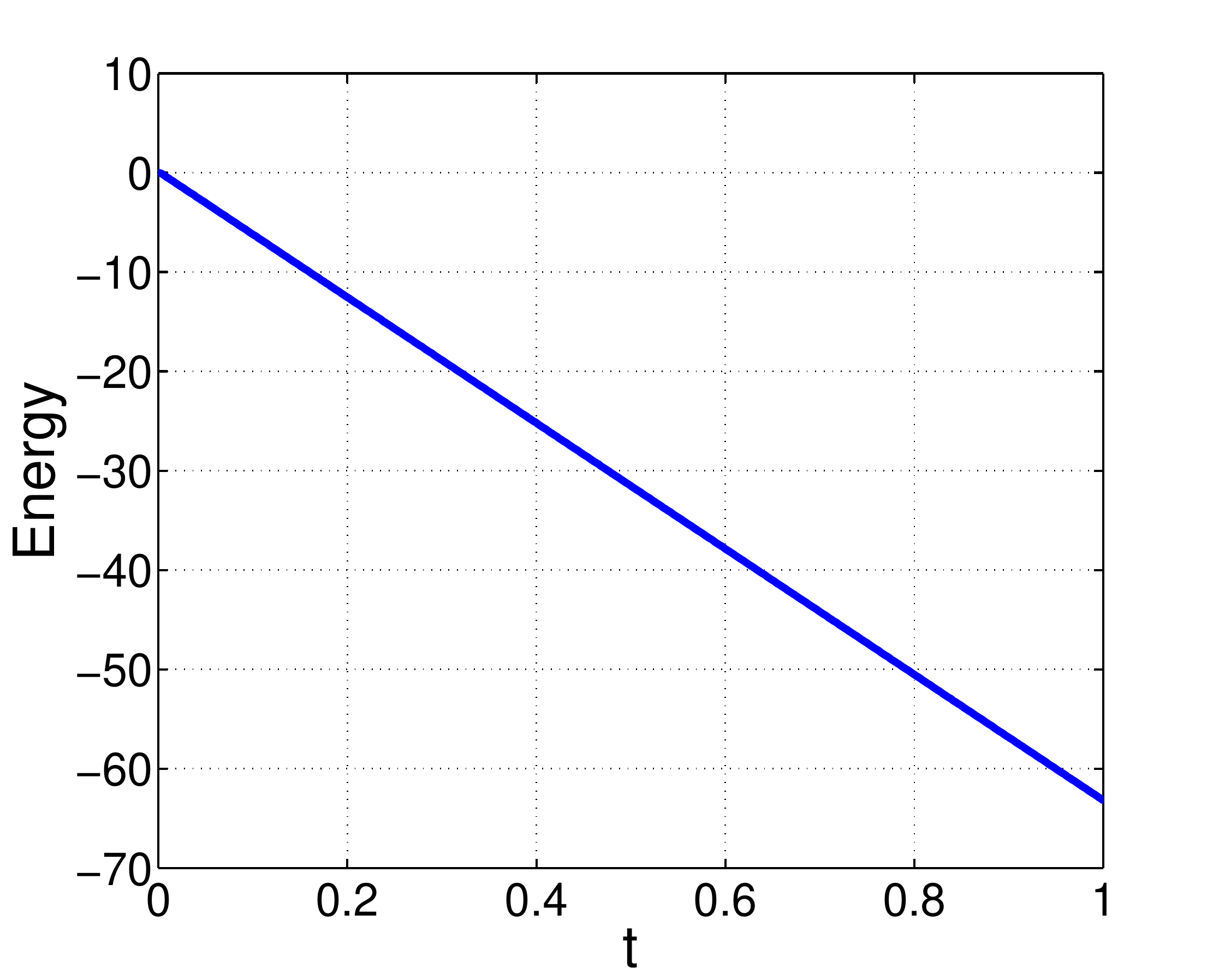}}
\subfloat[(f)]{\includegraphics[width=0.37\textwidth]{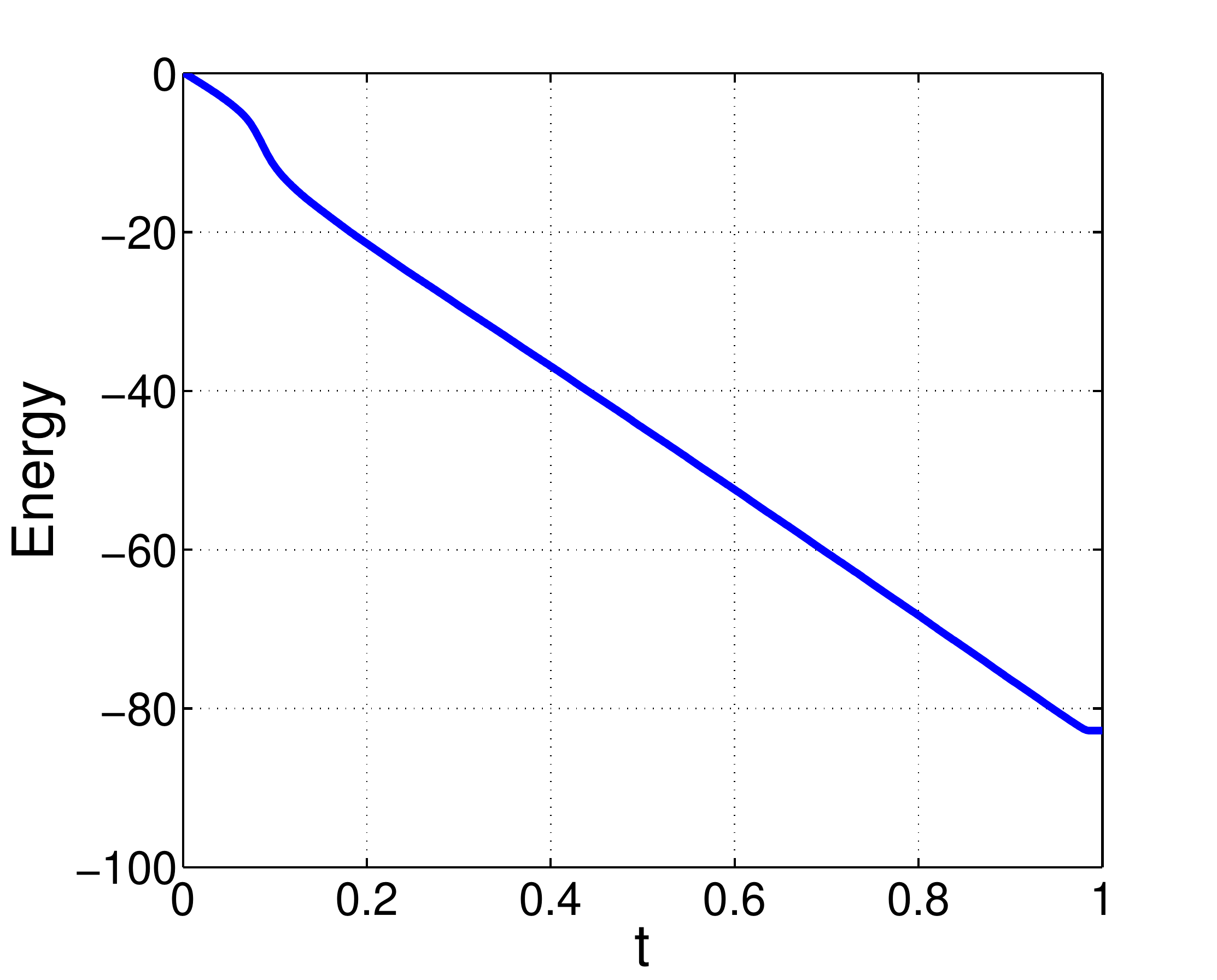}}
\subfloat[(i)]{\includegraphics[width=0.37\textwidth]{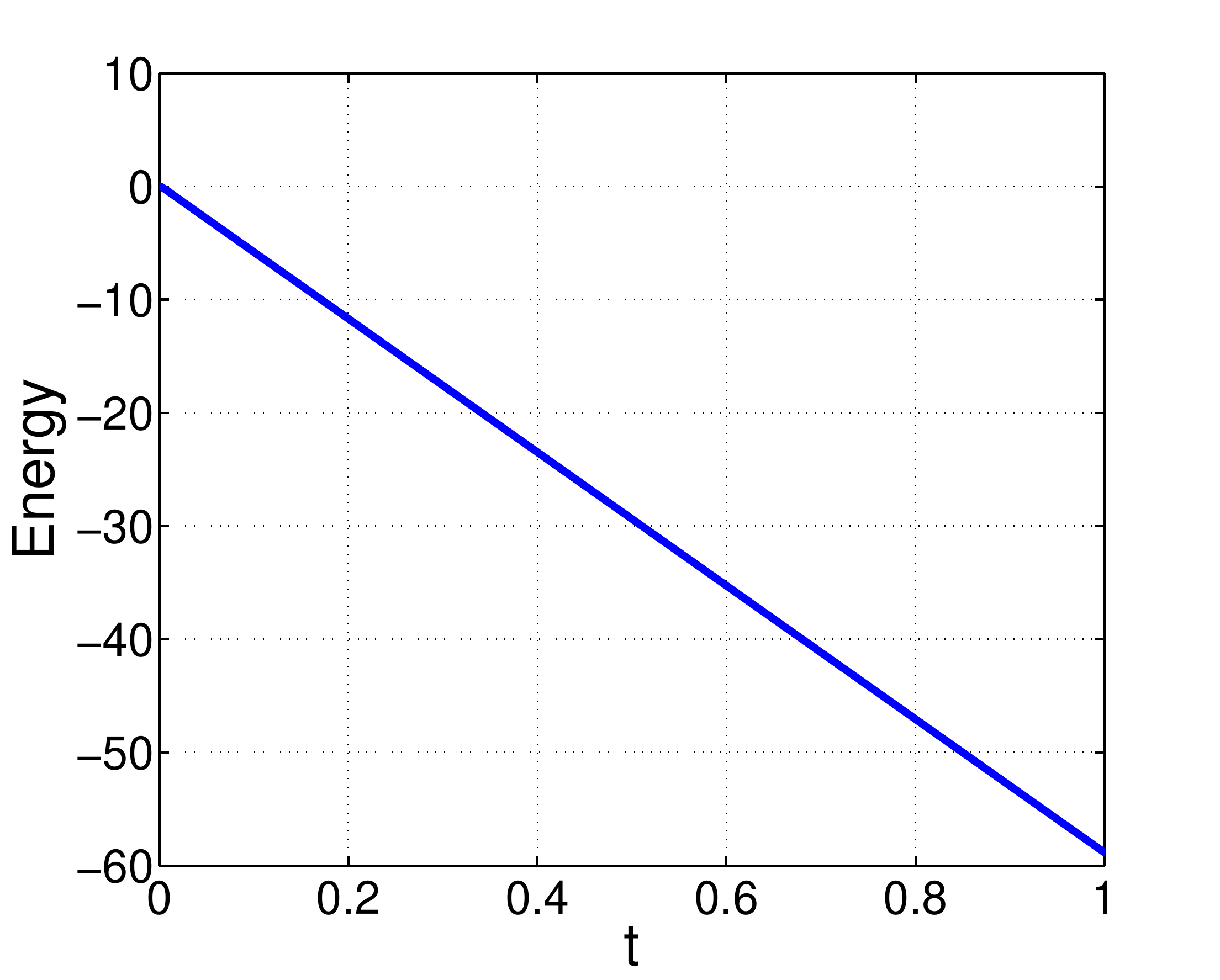}}
\caption{Schl\"{o}gl equation: Traveling wave solutions and energy plots; (a)-(c) on a uniform mesh with $N_I=120$ elements using quadratic dG basis functions, (d)-(f) on a moving mesh with $N_I=40$ elements using linear dG basis functions, (g)-(i) on a moving mesh with $N_I=40$ elements using quadratic dG basis functions.\label{ex3:plots}}
\end{figure}

Figure~\ref{ex3:mesh} shows that the steep wave fronts are captured well by the adaptive moving mesh method with linear and quadratic dG basis functions. On the other hand, Figure~\ref{ex3:plots} shows that the numerical solutions with quadratic dG basis functions give the correct wave speed, whereas for the linear case the numerical solutions move faster than the exact solutions. In Table~\ref{ex3:table}, we list the $L^2$-errors between the numerical and exact solutions at different time instances for different types of monitor functions. The correct wave speed is not captured by linear dG basis functions. Therefore  the corresponding $L^2$-errors are larger than the quadratic case.  Due the discontinuous nature of the dG discretization,  higher order basis functions can be implemented in a more flexible way in contrast to continuous finite elements, providing continuity requirement between the inter-element boundaries. The free energy of Schl\"ogl equation \eqref{energy}  decreases monotonically in time. The energy decreasing property is also captured numerically on uniform and moving meshes, as shown in Figure~\ref{ex3:plots} for the conditionally energy stable backward Euler method \cite{Hairer13edi}.

\begin{table}[htb!]
\centering
\caption{Schl\"{o}gl equation: $L^2$-errors with an adaptive moving mesh with $N_I=40$ elements.\label{ex3:table}}
\begin{tabular}{ l | c | l l l l l l}
Monitor & degree & t=0.001  & t=0.01 & t=0.25 & t=0.5 & t=0.75 & t=1 \\
\hline
Optimal  & 1 & 3.6e-03  & 1.6e-03 & 3.5e-01 & 4.3e-01 & 4.9e-01 & 5.3e-01  \\
Arc-length  & 1 & 9.3e-03  & 2.1e-02 & 4.9e-01 & 5.1e-01 & 5.4e-01 & 5.3e-01 \\
Curvature & 1 & 3.9e-03  & 5.2e-03 & 1.2e-01 & 1.8e-01 & 2.3e-01 & 2.5e-01 \\
\hline
Optimal  & 2 & 2.2e-04  & 3.3e-04 & 2.9e-03 & 4.4e-03 & 5.6e-03 & 6.7e-03 \\
Arc-length  & 2 & 1.1e-03  & 1.5e-03 & 2.2e-03 & 3.8e-03 & 8.8e-03 & 1.4e-03 \\
Curvature & 2 &  3.4e-04  &  2.3e-04 & 2.9e-03 & 5.2e-03 & 7.4e-03 & 9.5e-03 \\
\hline
\end{tabular}
\end{table}

\section{Conclusions}
\label{sec:conc}
In this paper we have developed an adaptive  discontinuous Galerkin moving mesh method for a class of one dimensional nonlinear PDEs. The moving mesh equations are solved using the static rezoning approach with the Lagrange dG basis functions in the Algorithm~\ref{alg}. Numerical results for problems with different nature of traveling waves demonstrate the accuracy and effectiveness of the moving mesh dG method. As a future study we aim to extend the dG moving mesh to two dimensional problems in space.

\end{document}